\documentclass[12pt]{amsart}

\pdfoutput=1

\usepackage{a4wide}
\usepackage{amsfonts,amsmath,amsthm, amssymb}
\usepackage[alphabetic]{amsrefs}
\usepackage{graphicx}

\newcommand{\pfi}[2]{\ensuremath{\partial_{#1}f_{#2}}}

\theoremstyle{plain}
\newtheorem*{main theorem}{Main Theorem}
\newtheorem*{theorem*}{Theorem}
\newtheorem*{emptytheorem*}{}
\newtheorem{theorem}{Theorem}
\newtheorem{proposition}{Proposition}
\newtheorem*{proposition*}{Proposition}
\newtheorem{lemma}{Lemma}

\newtheorem{corollary}{Corollary}
\newtheorem{claim}{Claim}

\theoremstyle{definition}
\newtheorem{definition}{Definition}
\newtheorem*{definition*}{Definition}

\theoremstyle{remark}

\newtheorem*{remark*}{Remark}

\begin{document}
\title[Boundary of hyperbolicity 
for H\'enon-like families]{The boundary of hyperbolicity \\ 
for H\'enon-like families}

\author{Yongluo Cao}
\address{Department of Mathematics, Suzhou University, 
Suzhou 215006, Jiangsu, P.R. China}
\email{{ylcao@suda.edu.cn}, 
{yongluocao@yahoo.com}}

\author{Stefano Luzzatto}
\address{Dept. of Mathematics, Imperial College, 180 Queen's 
Gate, London SW7 2AZ, UK}
\email{Stefano.Luzzatto@imperial.ac.uk}
\urladdr{http://www.ma.ic.ac.uk/\textasciitilde 
luzzatto}

\author{Isabel Rios}
\address{Universidade Federal Fluminense, Niteroi, RJ, Brazil.}
\email{rios@mat.uff.br}

\thanks{
IR was partially supported by CAPES, FAPERJ (Brazil) and EPSRC UK. 
YC was partially supported by  NSF(10571130) and NCET of China,
the Royal Society and EPSRC UK. SL was partially supported by  
NSF(10571130) and NCET of China.
This work was carried out at Imperial College London and
Suzhou University and we acknowledge the hospitality and
support of these institutions. We would also like to thank the referee 
for a careful reading of the paper and several very useful
suggestions which have improved the accuracy and 
presentation of the arguments.}   

\date{February 9, 2005, revised May 10, 2007}


\begin{abstract} 
 We consider \( C^{2} \) 
 H\'enon-like families of diffeomorphisms of \( \mathbb R^{2} \) and
 study the boundary of the region of parameter values for which the
 nonwandering set is uniformly hyperbolic. Assuming sufficient
 dissipativity, we show that the loss of hyperbolicity is caused by a 
 first homoclinic or heteroclinic tangency and that uniform 
 hyperbolicity estimates
 hold \emph{uniformly in the parameter} up to this bifurcation
 parameter and even, to some extent, at the bifurcation parameter. 
\end{abstract}    

\maketitle

\section{Introduction and statement of results} 
Our aim in this paper is to study the \emph{boundary of 
hyperbolicity} of certain families of two dimensional maps.

\subsection{H\'enon-like families}
We say that a family of \( 
    C^{2} \) plane 
    diffeomorphisms is called a \emph{
    H\'enon-like family} if 
    it can be written in the form 
    \[ 
   f_{a, b , \eta}(x,y) = (1-ax^{2}+y, bx) + \varphi(x,y,a) 
    \]
  where \( a\in\mathbb R \), \( b \neq 0 \) 
  and \( \varphi (x,y,a) \) is a  \( C^{2} \) ``perturbation''
  of the standard  \emph{H\'enon family} 
  \( h_{a,b}(x,y) = (1-ax^{2}+y, bx)  \) \cite{Hen76}
  satisfying 
  \[
  \|\varphi \|_{C^{2}(x,y,a)}\leq \eta.  
  \] 
In this paper we consider \( |b| \neq 0, \eta > 0 \) fixed
sufficiently small and investigate the dynamics as the 
parameter \( a \) is varied. For simplicity we shall therefore omit \( 
b \) and \( \eta \) from the notation and denote a H\'enon-like
family by \( \{f_{a} \} \). 
For future reference we remark that the inverse of \(
f_{a} \) is given by an equation of a similar form: 
\[ 
f^{-1}_{a} (x,y) = (y/b, x-1 + ay^2/b^2) + \tilde\varphi (x,y,a)
\]
where 
\( \|\tilde \varphi \|_{C^{2}(x,y,a)} \to 0\) as  
\( \|\varphi \|_{C^{2}(x,y,a)} \to 0 \). 
We shall suppose without loss of generality that 
 \[
 \|\tilde \varphi \|_{C^{2}(x,y,a)} \leq \eta.
 \]

\subsection{The boundary of hyperbolicity}

\subsubsection{Basic background}
H\'enon  and H\'enon-like families have been extensively
studied over the last almost 30 years. 
One of the earliest 
rigorous results on the subject is \cite{DevNit79} in which it was 
shown that the non-wandering set \( \Omega_{a,b} \) 
is uniformly hyperbolic for all \( b\geq 0 \) and 
all sufficiently large \( a \) (depending on \( b \)). 
On the other hand, for small \( b\neq 0 \) 
and \( a\lesssim 2 \) there exists positive probability of 
``strange attractors'' which contain tangencies between stable and 
unstable leaves. This was first proved in \cite{BenCar91}
for the H\'enon family and later generalized in \cite{MorVia93}
to H\'enon-like families, see also \cites{WanYou01, LuzVia03}. 
These attractors cannot be uniformly hyperbolic due to the presence of 
tangencies but turn out to satisfy weaker \emph{nonuniform} 
hyperbolicity conditions \cites{BenYou93,BenYou00, BenVia01}. 

\subsubsection{Complex methods}
More recently Bedford and Smillie  have described the transition
between these two regimes  for H\'enon families by identifying and 
describing some of the properties of the 
\emph{boundary of uniform hyperbolicity} \cite{BedSmi06}.   
In particular they show that for small 
\( |b| \), the nonwandering set is uniformly
hyperbolic up until the first parameter \( a \) at which a tangency
occurs between certain stable and unstable manifolds. 
Combining this with the statements contained in
\cite{BedSmi02} their results also
imply uniform bounds on the Lyapunov exponents of all invariant
probability measures at the bifurcation parameter \cite{Bed05}. 
Their methods
rely crucially on previous work \cite{BedSmi04} which in turn is
based on the polynomial nature of the H\'enon family, a feature
which allows them to consider its complexification and to apply
original and highly sophisticated arguments of holomorphic dynamics. 

\subsubsection{Real methods}
In this paper we develop a new and completely different strategy to the
problem, based purely on geometric ``real'' arguments, 
which have the advantage of applying to  general
\( C^{2} \) \emph{H\'enon-like} families. We also obtain 
the analogous \emph{uniformity} results by showing that
the hyperbolicity expansion and contraction rates are uniform right up
to the point of tangency and that even \emph{at} the point of tangency
some strong version of nonuniform hyperbolicity continues to hold: all
Lyapunov exponents of all invariant measures are uniformly bounded
away from 0.

\begin{theorem}
For all \( |b| >  0 \) and \( \eta >0  \) sufficiently small we have
the following property. For every H\'enon-like family \( \{f_{a}\}_{\
a\in\mathbb R} \) of
plane diffeomorphisms there exists a unique \( a^{*} \) such that 
\begin{enumerate} 
    \item For all \( a> a^{*} \) the nonwandering set \( \Omega_{a} \)
    is uniformly hyperbolic;
    \item For \( a=a^{*} \) the nonwandering set \( \Omega_{a^{*}} \) 
    contains an orbit of tangency but is ``almost uniformly hyperbolic''
    in the sense that all Lyapunov exponents of all invariant
    probability measures supported on \( \Omega \) are uniformly
    bounded away from 0. 
\end{enumerate}
Moreover, the bounds on the expansion and contraction rates for all \( 
a\geq a^{*} \) are independent of \( a \) and of the family. 
\end{theorem}

\subsubsection{Singular perturbations}
We remark that this is not the only existing definition of 
H\'enon-like  in the literature. One standard approach is to
consider  ``singular'' perturbations of the limiting one-dimensional
map corresponding to the case \( b=0 \):
\[ 
f_{a}(x,y) = (1-ax^{2}, 0) + \varphi_{a}(x,y).
\]
This formulation however has some slight technical issues.  For
example,  one 
cannot assume that \( \|\varphi_{a}\|_{C^{2}} \) is small on all of \( \mathbb 
R^{2} \) since that
would violate the requirement that \( f_{a} \) be a global
diffeomorphism of \( \mathbb R^{2} \). This can be dealt with by
restricting our attention to some compact region, say \( [-2, 2]
\times [-2,2] \), and supposing only that \(
\|\varphi_{a}\|_{C^{2}}\leq \eta \) in this region. Our arguments apply
in this case also and yield a more local result on the hyperbolicity of the
nonwandering set restricted to \( [-2, 2] \times [-2,2] \). 

\subsection{Basic definitions}\label{basicdefs}

\subsubsection{Nonwandering set}
We recall that a point \( z \) belongs to the \emph{nonwandering} set \(
\Omega \) of \( f \) if it has the property that for every
neighbourhood \( \mathcal U \) of \( z \) there exists some \( n\geq 1 \)
such that \( f^{n}(\mathcal U)\cap \mathcal U \neq \emptyset \). The
nonwandering set is always invariant and closed (and thus if bounded, 
also compact). 

\subsubsection{Uniform hyperbolicity}
We say that a compact invariant set \( \Omega \) is \emph{uniformly hyperbolic} 
(with respect to \( f \)) if 
there exists constants \( C^{u}, C^{s}>0, \lambda^{u}>0>\lambda^{s} \) 
and a \emph{continuous} 
decomposition \( T\Omega=E^{s}\oplus E^{u} \) of the tangent bundle 
such that for every \( x\in \Omega \), every 
non-zero vector \( v^{s}\in E^{s}_{z} \) and \( v^{u}\in E^{u}_{z} \) and every 
\( n\geq 1 \) we have 
\begin{equation}\label{UH}
\|Df^{n}_{z}(v^{s})\|\leq C^{s}e^{\lambda^{s} n} \quad\text{and} \quad 
\|Df^{n}_{z}(v^{u})\|\geq C^{u}e^{\lambda^{u} n}. 
\end{equation}
By standard hyperbolic theory, the stable and unstable subspaces \( 
E^{s}_{z}, E^{u}_{z} \) are tangent everywhere to the stable and 
unstable manifolds. In particular uniform hyperbolicity is 
incompatible with the presence of any tangencies in \( \Omega \) 
between any stable 
and any unstable invariant manifolds associated to points of \( \Omega 
\). 

\subsubsection{Nonuniform hyperbolicity}
A weaker notion of hyperbolicity can be formulated in terms of
invariant measures. For simplicity we restrict our discussion to the
two-dimensional setting, as relevant to the situation we consider in
this paper.
Let \( \mu \) be an \( f \)-invariant ergodic probability measure 
    with support in some compact invariant set \( \Omega \). 
By Oseledec's Ergodic Theorem \cite{Ose68} there exist constants 
\( \lambda^{u} \geq \lambda^{s} \) and a 
measurable decomposition  
\( T\Omega=E^{s}\oplus E^{u} \) such that for \( \mu \)-almost every \( 
z \) and  every non-zero
vector \( v^{s}\in E^{s}_{z} \) and \( v^{u}\in E^{u}_{z} \) 
we have 
\begin{equation}\label{NUH}
\lim_{n\to\infty} \frac{1}{n}\log\|Df^{n}_{z}(v^{s})\| = \lambda^{s} 
\quad\text{and}\quad
\lim_{n\to\infty} \frac{1}{n}\log\|Df^{n}_{z}(v^{u})\| = \lambda^{u}.
\end{equation}
The constants \( \lambda^{s} \) and \( \lambda^{u} \) are called the
\emph{Lyapunov exponents} associated to the measure \( \mu \).
We say that \( \mu \) is \emph{hyperbolic} \cites{Pes76,Pes77} if 
\[
\lambda^{u}>0>\lambda^{s}.  \]
Clearly \eqref{UH} implies \eqref{NUH} for any \( \mu \). 
The converse however is false in 
general: the measurable decomposition may not extend to a continuous 
one on all of \( \Omega \) and the exponential expansion and 
contraction in \eqref{NUH} implies only a limited version of \eqref{UH} in 
which the constants \( C^{s}, C^{u} \) are measurable functions 
of  \( x \) and not  bounded away from 0. 
This definition of hyperbolicity in terms of Lyapunov exponents is
sometimes called  \emph{nonuniform 
hyperbolicity} and is consistent in principle with the existence 
of tangencies between stable and unstable manifolds.  

\subsubsection{The boundary between uniform and nonuniform
hyperbolicity}
In general there may be many ergodic 
invariant probability measures supported in \( \Omega \) of which some
may be hyperbolic and some not. Even if they are all 
hyperbolic the corresponding Lyapunov exponents may not be uniformly 
bounded way from \( 0 \). The situation in which all Lyapunov 
exponents of all ergodic invariant measures are uniformly bounded away 
from zero is, in some sense, as ``uniformly hyperbolic'' as one can 
get while admitting the existence of tangencies. This situation 
can indeed occur, for example in the present context of H\'enon-like
maps.  A first example of a 
set satisfying this property was given in \cite{CaoLuzRioTan}.

\subsection{A one-dimensional version}

After completing the proof of the Theorem 1 we realized that 
\emph{much simpler} versions of our arguments yield an analogous, new 
and non-trivial, result in the context of one-dimensional maps.  We
explain and give a precise formulation of this result. We consider
first the quadratic family 
\[ 
h_{a}(x) = 1-ax^{2}.
\]
We choose this particular parametrization for convenience and
consistency with our two dimensional results, but any choice of smooth
family of  unimodal or even multimodal maps with negative Schwarzian
derivative would work in exactly the same way. It is well known that
for \( a>2 \) the nonwandering set \( \Omega_{a} \) is uniformly
expanding although we emphasize here that this depends crucially on
the negative Schwarzian derivative property. The negative Schwarzian
property is not robust with respect to \( C^{2} \) 
perturbations and standard methods do not therefore yield this
statement for such perturbations. 

\begin{theorem}
    There exists a constant \( \eta>0 \) such that if a family \(
    \{g_{a} \}\) of \(  C^{2} \) one-dimensional maps satisfies 
    \[ 
    \|g_{a}-h_{a}\|_{C^{2}}\leq \eta
    \]
    then there exists a unique parameter value \( a^{*} \) such that
    \begin{enumerate}
\item For all \( a>a^{*} \) the non-wandering set \( \Omega \) is
uniformly hyperbolic; 
\item For \( a=a^{*} \) the Lyapunov exponents of all ergodic
invariant probability measures are all positive and uniformly bounded 
away from 0.
\end{enumerate}
Moreover the rates of expansion and the bound on the Lyapunov
exponents are uniform, independent of the family and of the parameter.
\end{theorem}

The proof of this result is  exactly the same as that of
Theorem 1 but hugely simpler as all more geometrical arguments
concerning curvature etc become essentially trivial. 

We emphasize that the
uniform expansivity of \( \Omega_{a} \) for a particular parameter
value \( a>2 \) is of course robust under
sufficiently small perturbations of \( f_{a} \), by standard
hyperbolic theory. However this approach \emph{requires the size of the
perturbation to depend on the parameter }\( a \) and in particular to
shrink to zero as \( a \) tends to \( 2 \).  The crucial point of our 
approach, both in this one-dimensional setting, as in the
two-dimensional setting is that the size of the perturbation
\emph{does not} depend on the parameter.

\subsection{Overview of the paper}

We have divided our argument into three main sections. 
In 
Section \ref{proofnonwan}
we analyze the geometric structure of
stable and unstable manifolds of the two fixed points and define the
parameter \( a^{*} \) as the first parameter for which a tangency
occurs between some compact parts of these manifolds. We also identify
a region \( \mathcal D \) which we show contains the non-wandering set.
In Section \ref{sectionhyp} 
we define a 
``critical neighbourhood'' \( \Delta_{\varepsilon} \) outside of which
our maps are uniformly hyperbolic by simple perturbation arguments. 
However \(\Delta_{\varepsilon} \) does contain points of \( \Omega \) 
and thus we cannot ignore this region. 
To control the hyperbolicity in \( \Delta_{\varepsilon} \) we
 introduce the notions of Hyperbolic Coordinates and Critical Points 
which form the key technical core of our approach. 
Finally, in Section 
\ref{sectionhypest} we apply these techniques to prove the required
hyperbolicity properties.

\section{The non-wandering set}\label{proofnonwan}

In this section 
we define the parameter \(a^{*}\) as in the statement of our main
Theorem, 
and show that for \( a\geq a^{*} \) the
nonwandering set is contained in the closure if the unstable manifold 
of a hyperbolic fixed point restricted to a certain compact region of \( 
\mathbb R^{2} \). 

\subsection{The parameter \protect\( a^{*} \protect\)}

We define the bifurcation parameter \( a^{*} \) below as the first
parameter for which there is a tangency between certain compact parts 
of the stable and unstable manifolds of the fixed points. This does
not immediately imply that it is a first parameter of tangency though 
this will follow from our proof of the fact that the nonwandering set 
is uniformly hyperbolic for all \( a> a^{*} \). 

\subsubsection{Fixed points and invariant manifolds 
for the one-dimensional limit}
For the
endomorphisms \(h_{a}=h_{a,0}\) with \(a\geq 2\), 
there are two fixed points,
\[
p_a = \frac{-1+\sqrt{1+4a}}{2a}> q_a = \frac{-1-\sqrt{1+4a}}{2a}
\] 
both hyperbolic. 
For the special parameter value \( a=2 \), to simplify the notation
below, we write 
\[ f_{*}=h_{2,0}, \text{ and denote the two fixed points by
} 
p_{*}=(1/2, 0)  \text{ and }  q_{*}=(-1,0). 
\]
Since \( q^{*} \) and \( p^{*} \) are repelling in the horizontal
direction, their stable sets are simply their preimages: 
\[ W^{s}(q^{*}) = \bigcup_{n\geq 0}f_{*}^{-n}(q^{*}) 
\quad\text{and} \quad 
W^{s}(p^{*})=\bigcup_{n\geq 0}f_{*}^{-n}(p^{*}). 
\]
In particular these sets contain the following curves
\[ 
f_{*}^{-1}(q_{*}) = \{(x, y):
f_{*}( x,  y) = (1-2 x^{2}+ y, 0) = (-1, 0)\}
= \{ y = 2 x^{2} - 2\} 
\]
and
\[
f_{*}^{-2}(q_{*}) = \{(x, y):
f_{*}( x,  y) = (1-2 x^{2}+ y, 0) = (1, 0)\}
= \{ y = 2 x^{2} \} 
\]
\begin{figure}
    \includegraphics[width=0.3\textwidth]{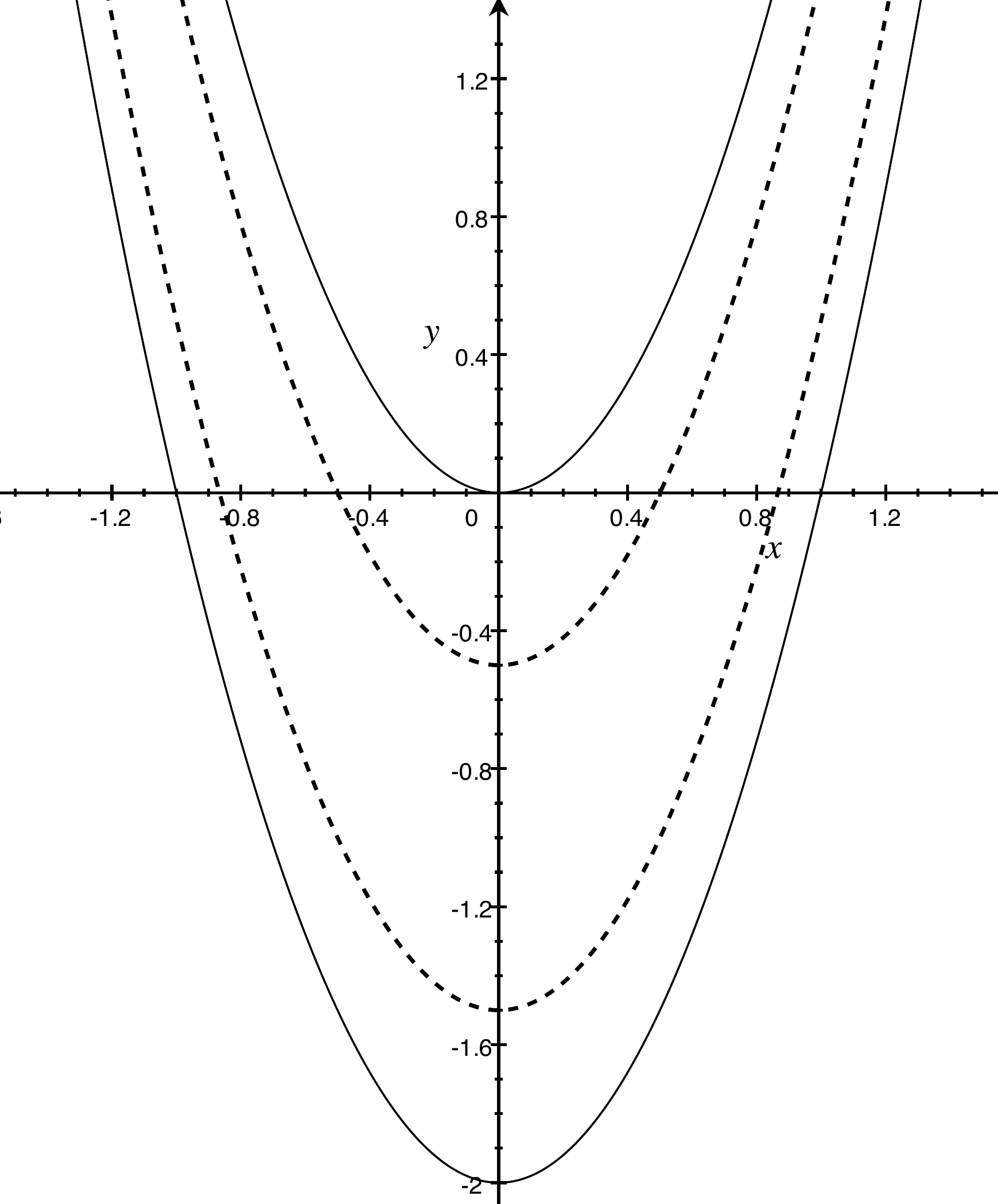}
   \caption{First two ``generations'' of \( W^{s}(q_{*}) \)
  and  \( W^{s}(p_{*}) \).}
  \label{Wsb=0}
  \end{figure}
The first preimage of \( q^{*} \) 
is a parabola with a minimum at \( (0, -2) \) and 
intersecting the \( x \)-axis at  \( x=\pm 1 \), and having
slope equal to \( -4 \) at the point \( q_{*}=(-1,0) \) and the second
is a parabola with a minimum at \( (0,0) \).
Similarly we can compute 
\[ 
f_{*}^{-1}(p_{*}) = \{ z=( x,  y): 
f_{*}( z) = (1-2 x^{2}+ y, 0) = (1/2, 0)\}
= \{ y = 2 x^{2} - 1/2\}
\]
which is a parabola with a minimum at \( (0, -1/2) \), 
intersecting the \( x \)-axis at \( x= \pm 1/2 \) and  having 
slope equal to \( 2 \) at the point \( p_{*}= (1/2, 0) \), and 
\[ 
f_{*}^{-2}(p_{*}) = \{ z=( x,  y): 
f_{*}( z) = (1-2 x^{2}+ y, 0) = (-1/2, 0)\}
= \{ y = 2 x^{2} - 3/2\}
\]
which is a parabola with a minimum at \( (0, -3/2) \).

The unstable manifolds going \( W^{u}(q_{*}) \) and \( W^{u}(p_{*}) \)
can be defined and computed in a similar way and are easily seen to be
horizontal.

\subsubsection{Fixed points for H\'enon-like families}
Consider first the \emph{H\'enon family} 
\(h_{a,b}(x,y)=(1-ax^{2}+y,bx).\)
For \(b\neq 0\), \(h_{a,b}\) is a diffeomorphism. 
The hyperbolicity of the fixed points  implies that there
exists a neighbourhood of the set \( \{ (a,0):a\geq 2 \}\) 
corresponding to
pairs of parameters for which there is an
\emph{analytic continuation} \( q_{a,b}, p_{a,b} \)  
as hyperbolic fixed points of \(h_{a,b}\). Considering \(\eta\) small,
we also have that the analytic continuations \(q_{f_{a}}\) and
\(p_{f_{a}}\) are also well defined and hyperbolic. For simplicity we 
shall often just refer to these two points as \( q, p \) leaving
implicit their dependence on \( f \). 

Explicit formulas for \( q_{a,b}, p_{a,b} \)  
 can be easily derived from the equation 
\(
( 1-ax^{2}+y, bx )=   (x,  y ) 
\)
but these would not be particular useful. Instead we just observe that
the fixed points must lie on the line 
\( 
\{y=bx\}
\)
and so in particular this means that for \( a\approx 2 \) and \( b
\gtrapprox 0 \)  the vertical coordinates 
of \( q_{a,b} \) and \( p_{a,b} \) are negative and positive
respectively, and the converse for \( b \lessapprox 0 \). Clearly the 
same holds for \(q=q_{f_{a}}\) and \(p=p_{f_{a}}\) if \( \eta \) is
sufficiently small.  
Moreover, the determinant of \( h_{a,b} \) is given by 
\[ 
\det Dh_{a,b}=\det \begin{pmatrix} 
-2ax & 1 \\ b & 0 \end{pmatrix} = -b.
\]
In particular the determinant is constant and negative if \( b \) is
positive and positive if \( b \) is negative. We thus refer to the
case \( b>0 \) as the \emph{orientation-reversing} case, and the case \( 
b<0 \) as the \emph{orientation-preserving} case.
\begin{figure}[h]
  \includegraphics[width=0.45\textwidth]{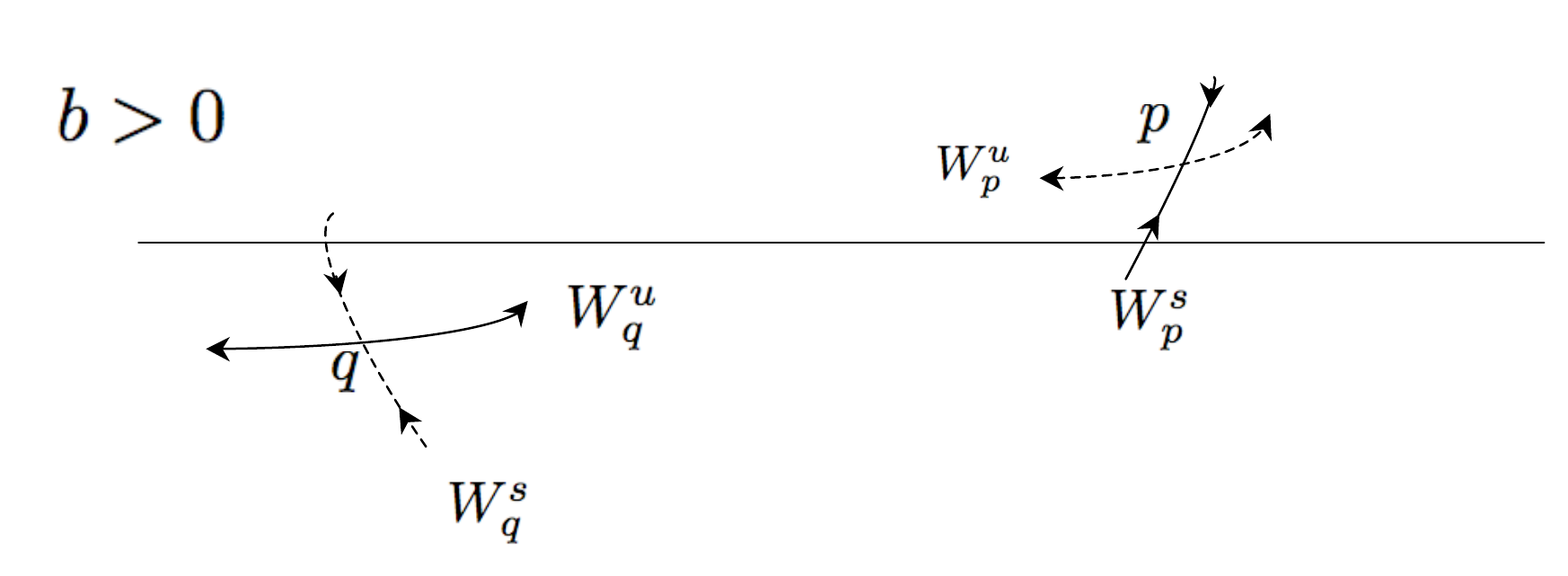}
  \includegraphics[width=0.45\textwidth]{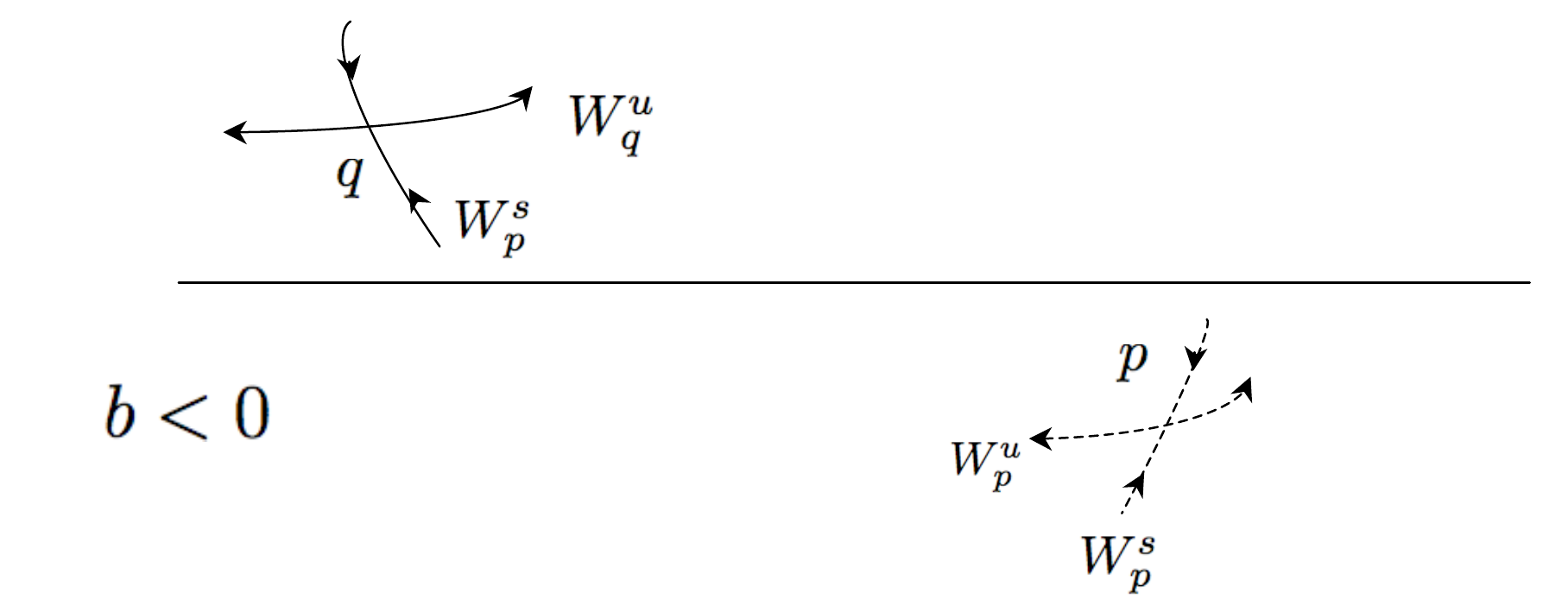}
\caption{Fixed points and their local stable and unstable manifolds
for the orientation-reversing (\( b>0 \)) and the
orientation-preserving (\( b<0 \)) case
(dotted curves indicate negative eigenvalues)}
\label{fixedpoints}
\end{figure} 
Recall that the determinant of a matrix is the product of the
eigenvalues, and thus in particular, the sign of the determinant 
has implications for the sign of the eigenvalues which, as we shall
see,  in turn has implications for the geometry of the stable
and unstable manifolds of the fixed points. For \( b=0 \) the fixed
points \( p_{*} \) and \( q_{*} \) have derivatives \( 4 \) and \( -2 \)
respectively, and thus, for \( b\neq 0 \) and \( \eta \) small, the
expanding eigenvalues of \( p \) and \( q \) are \( \approx 4 \) and \( 
\approx -2 \) respectively. This implies that for \( b \gtrapprox 0
\), the orientation-reversing case, the contracting eigenvalues of \( 
q \) and \( p \) must be \( <0 \) and \( >0 \) respectively, while for
\( b<0 \), the orientation preserving case, they must be \( >0 \) and \( 
<0 \) respectively. The two situations are illustrated in Figure
\ref{fixedpoints} with dashed lines showing the invariant manifolds
corresponding to negative eigenvalues.

\subsubsection{Analytic continuation of stable and unstable manifolds}
\begin{figure}[h]
    \includegraphics[width=9cm]{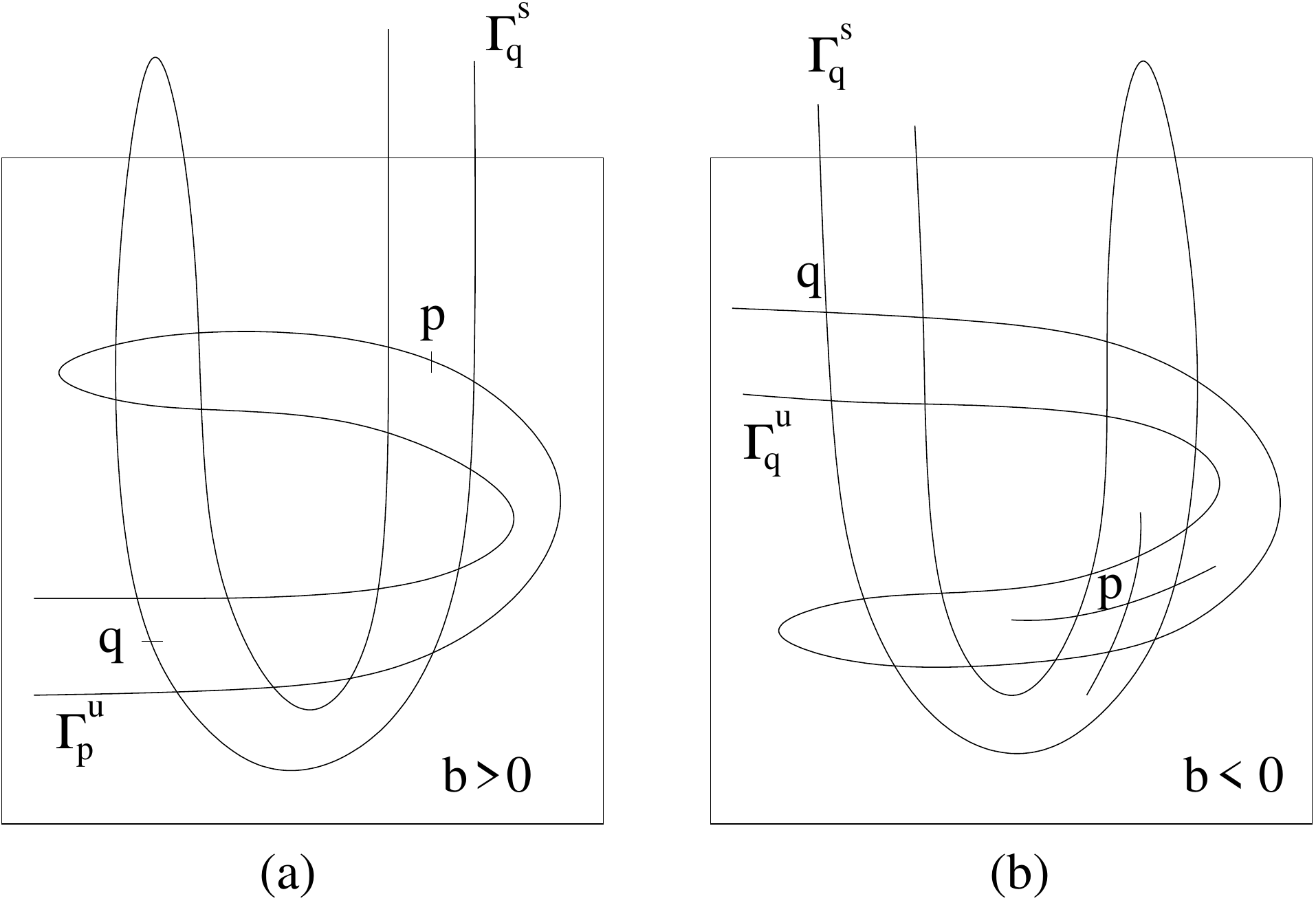}
  \caption{Invariant manifolds for \(a>a^{*}\)}
  \label{intersection}
\end{figure}

By classical hyperbolic theory, compact parts of the stable manifolds 
  depend continuously on  the map (see e.g.  \cite{PalMel82}).
  Therefore, for small \( b \) and small \( \eta \) 
  the analytic continuations \( q, p \)
  of the fixed points \( q_{*} \) and \( p_{*} \) 
  have stable and unstable manifolds which are close to those computed
  above for the limiting case. Elementary calculations show that the 
actual geometrical relations between these continuations depend on 
whether we consider the orientation reversing 
(\( b>0 \)) or the
orientation preserving (\( b<0 \)) case, and are as illustrated in
Figure \ref{intersection}.
We let 
\[
\Gamma^{u}_{a}(p)\subset W^{u}_{a}(p), \quad 
\Gamma^{s}_{a}(q)\subset
W^{u}_{a}(q),  \quad  
\Gamma^{u}_{a}(q)\subset W^{u}_{a}(q),
\] 
denote the compact parts of the stable and unstable manifold 
as shown in Figure \ref{intersection} and notice that,
in particular, since for \( b=0 \) and  \(a>2\)
the unstable manifold of \(p_{a}\) and \( q_{a} \) 
extend to the whole 
of the line, for each \( a>2 \) and \( b>0 \) sufficiently
small
we have that \(W^{u}_{loc}(p)\) crosses  \(W^{s}_{loc}(p)\) four
times, and for each \( a>2 \) and \( b<0 \) sufficiently small we have 
that \( W^{u}_{loc}(q) \) crosses \( W^{s}_{loc}(q) \) four times, and
also we can ensure that the compact parts defined above and in the
Figure intersect transversally. 
  Again this continues to hold
 also for a H\'enon-like family for sufficiently small \(
  \eta \). 
  
  \subsubsection{Definition of  \protect\( a^{*} \protect\)}\label{astar}
We are now ready to define the parameter \( a^{*} \).  We 
fix \( b\neq 0 \). \\
For an \emph{orientation-reversing }(\( b>0 \))
 H\'enon-like family \( f_{a} \), we let 
\[ 
a^{*}=\inf\{a: \Gamma^{s}_{a}(p) \text{ and } \Gamma^{u}_{a}(q) 
\text{ intersect transversally }\}.
\]
For an \emph{orientation-preserving} \( b<0 \) 
H\'enon-like family \( f_{a} \), we let 
\[ 
a^{*}=\inf\{a: \Gamma^{s}_{a}(q) \text{ and } \Gamma^{u}_{a}(q) 
\text{ intersect transversally }\}.
\]
We also  define a parameter 
\[ 
\hat a < a^{*} 
\]
as the \( \inf \)
  of parameters \( a \) for which 
\(W^{u}_{loc}(p)\) crosses  \(W^{s}_{loc}(p)\) 
  four  times (\( b>0 \)) or 
  \(W^{u}_{loc}(q)\) crosses  \(W^{s}_{loc}(q)\) 
    four  times (\( b<0 \)).  
  Clearly this is a weaker condition and thus \(  a^{*}\geq 
  \hat a\).  Notice that \( a^{*} \) and \( \hat a \) converge to \(
  a=2 \) as \( b \) and \( \eta  \) tend to 0.

  \subsection{Localization of the nonwandering set}

\begin{figure}[h]
      \includegraphics[width=4cm]{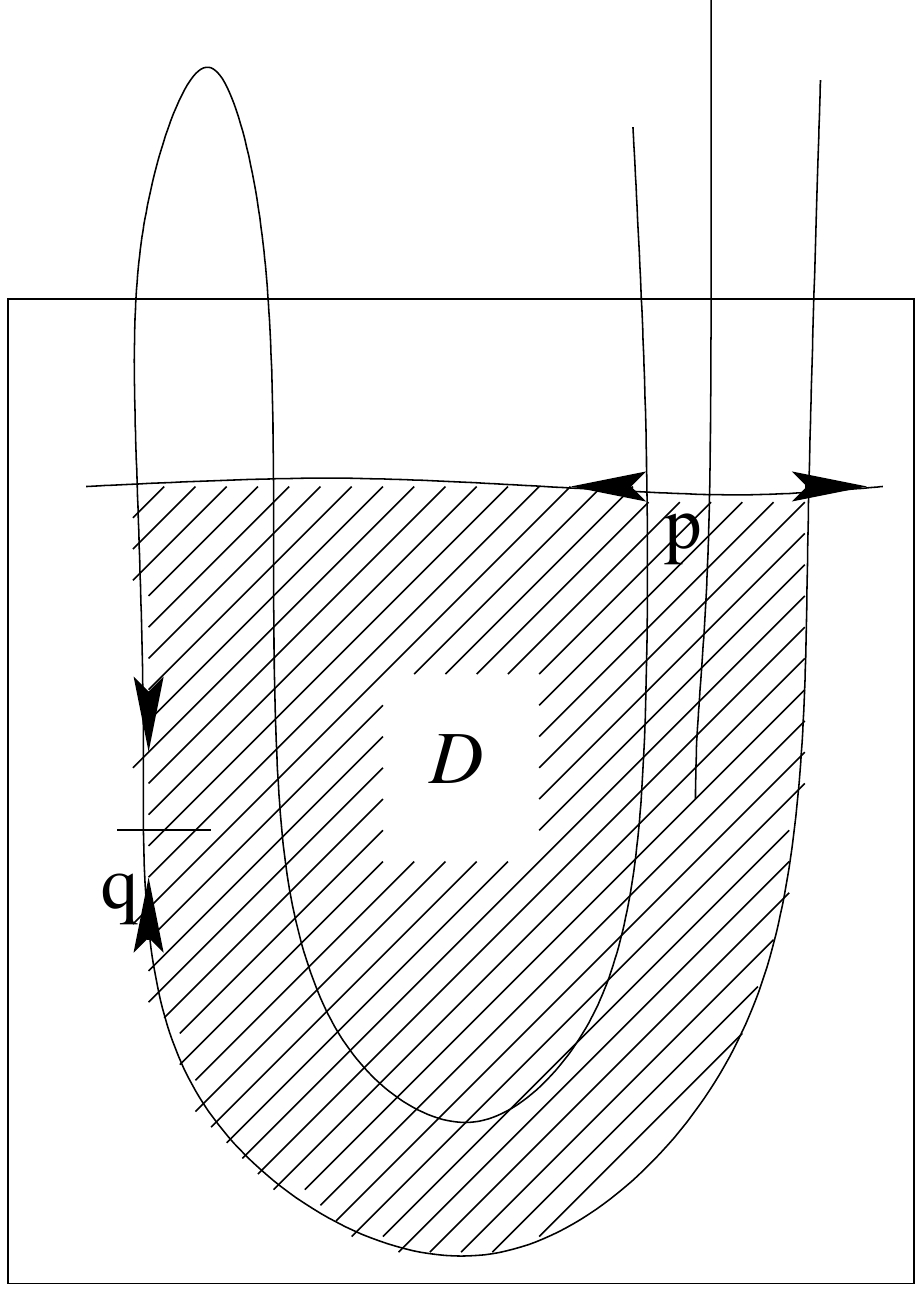}
    \caption{The region \( \mathcal D \)}
    \label{RegionD}
  \end{figure}

  In this section we carry out a detailed geometrical study aimed at 
  showing that the nonwandering set is contained in a relatively 
 restricted region. To prove hyperbolicity we will  then just have to 
 focus our efforts in this region.  For the moment we restrict
 ourselves to the orientation-reversing case. At the end we shall
 remark how the orientation-preserving case follows by identical
 arguments with a few minor changes of notation. 
First of all we let 
  \(
  \mathcal D 
  \) 
  denote the closed topological disc
  bounded by compact pieces of the \( W^{u}(p) \)
  and \( W^{s} (q)  \) as shown in Figure \ref{RegionD}.
The main result of this section is the following 
\begin{proposition}\label{nonwan}
For all \( a> \hat a \) we have 
\[
\Omega \subset \overline{W^{u}(p)} \cap  \mathcal D
\cap \{[-2,2] \times (-4b, 4b)\}.
\]
 \end{proposition} 
   
 \begin{remark*}
 In this paper we are interested in parameters \( a \geq  a^{*} (\geq  \hat a) \),
 but
 it is worth observing that if 
 follows from Proposition \ref{nonwan} that for all \( a\in (\hat
 a, a^{*}) \), and so in particular for a certain range of parameter
 values which  may contain multiple tangencies 
 the recurrent dynamics is  captured to some extent
 by the dynamics on \( W^{u}(p) \).
 This includes in particular
 all complex 
 dynamical phenomena associated to the unfolding of the tangency 
 at the parameter \( a^{*} \) (indeed, this includes the range of
 parameters considered by Benedicks-Carleson in \cite{BenCar91}. 
  \end{remark*}

\begin{figure}[h]
    \includegraphics[width=9cm]{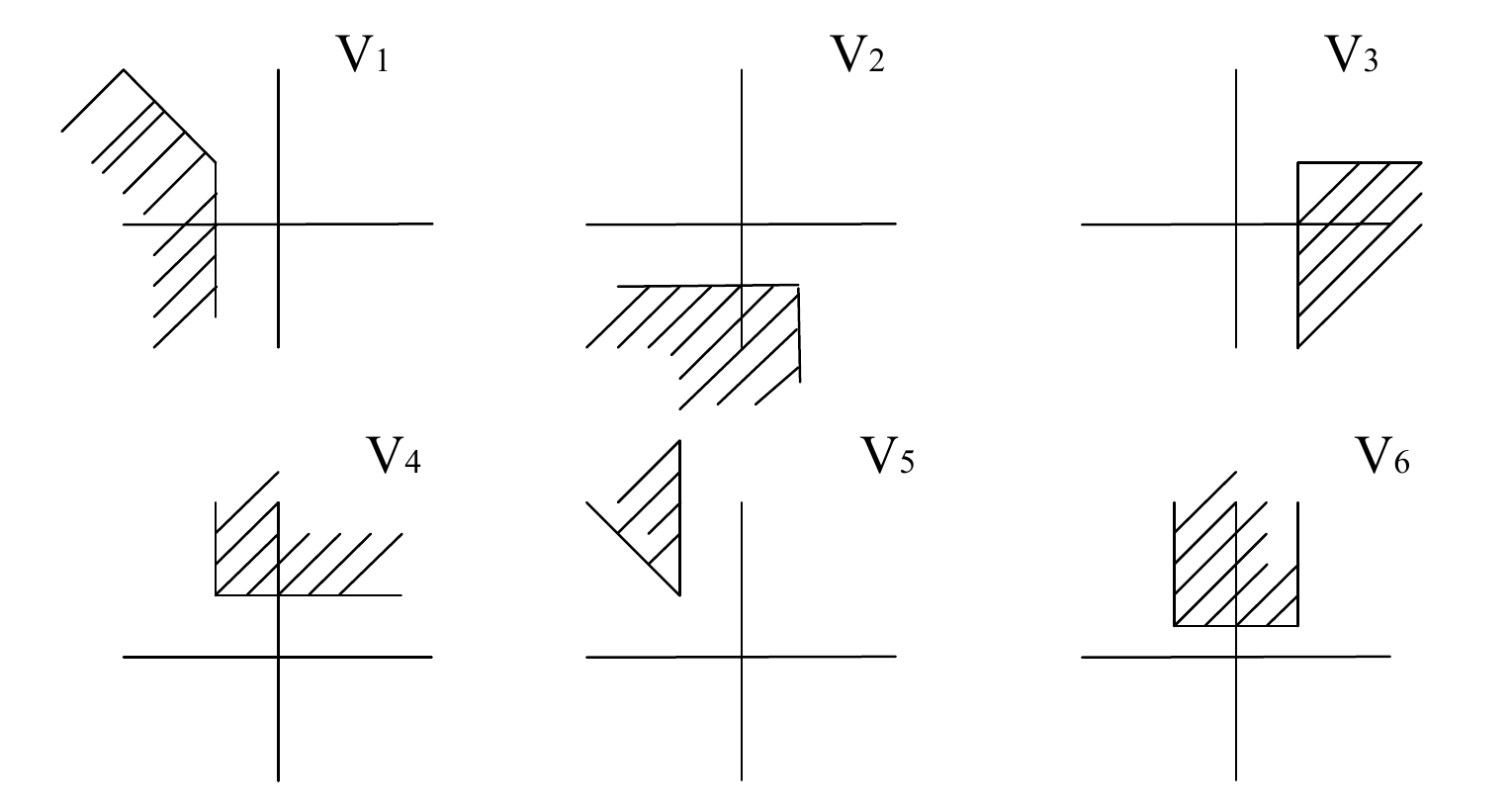}
  \caption{Regions \( V_{1} \) to \( V_{6} \)}\label{RegionV}
\end{figure}

We split the proof of Proposition \ref{nonwan} 
in  several Lemmas. Once again we deal first with the case \( b>0 \). 
At the end of the proof we indicate the necessary modifications in
order to deal with the case \( b< 0 \). 
We first define a relatively ``large'' region \( R \)
and show that  \( \Omega\subset R \) and then show in separate 
arguments that \(  \Omega \subset \mathcal D\) and 
\( \Omega\subset \overline{W^{u}(p)} \), and finally refine our
estimate to obtain the statement in the Proposition.
Let 
\[
R= (-2, 2)\times (-4, 4b) \subset \hat R = 
(-2, 2)\times (-4, 2) 
\]
We also define the following
 6 (overlapping)  regions (see Figure \ref{RegionV}):
\begin{align*}
       V_{1}&=\{(x,y): x\leq -2, y\leq |x|\}, \\
	V_{2}&=\{(x,y): x\leq 2, y\leq -4\},\\
	V_{3}&=\{(x,y): x\geq 2, y\leq 2\},\\
	V_{4}&=\{(x,y): x\geq -2, y\geq 2\},\\
	V_{5}&=\{(x,y): x\leq -2, y\geq |x|\},\\
       V_{6}&= \{(x,y): |x|\leq 2, y\geq 4b\}
\end{align*}
 
Then 
\[
\hat R= \mathbb R^{2}\setminus
(V_{1}\cup\dots\cup V_{5}) \text{ and }
R = \mathbb R^{2}\setminus
(V_{1}\cup\dots\cup V_{6})
\]
 
We prove the following two statements. 

 \begin{lemma}\label{1.1}
   \( \Omega \subset  R\).
    \end{lemma}
    \begin{proof}
We  show that the orbit of every 
point \( (x,y)\in V_{i} \), \( i=1,\ldots,6 \) is unbounded in either backward
or forward time.  This implies in particular that no such point is
nonwandering. 
For \( n \in \mathbb Z \), let \( (x_{n}, y_{n})=
f_{a}^n(x,y) \). We shall use repeatedly the fact that \( a\approx 2 \) 
and  \( b \approx 0 \). 

For \( (x,y) \in V_{1} \) we have \( x\leq -2 \) and therefore 
 \( x_{1}=1-ax^{2}+y +\varphi_{1}(x,y,a)\leq 1-ax^{2}+|x| +\eta \leq -2 \) and
 \( y_{1}=bx +\varphi_{2}(x,y,a) \leq -2b+\eta  <
|x_{1}| \), as long as \(\eta \) is sufficiently small.
Thus \( (x_1,y_1) \in V_{1} \), and \( |x_{1}| \geq ax^{2}-|x|-1 -\eta
\geq 2|x| \). 
Repeating the calculation we have
\( |x_{n}|\geq 2^{n}|x| \) and so \( |x_{n}|\to\infty \).

For \( (x,y)\in V_{2} \) we have \( x_{1}=1-ax^{2}+y
+\varphi_{1}(x,y,a)\leq -2\) and \(
y_{1}=bx +\varphi_{2}(x,y,a)\leq 2b +\eta < |x_{1}| \). Thus \( (x_{1},y_{1})\in V_{1} \) and so  
\( |x_{n}|\to\infty \).

Similarly, for \( (x,y)\in V_{3} \) we have
\( x_{1}=1-ax^{2}+y+\varphi_{1}(x,y,a) \leq 1- 2\cdot 2^{2} + 2+\eta \leq -2\) and \(
y_{1}=bx +\varphi_{2}(x,y,a)\leq 2b +\eta < |x_{1}| \). Thus \( (x_{1},y_{1})\in V_{1} \) and we
argue as above.

For \( (x,y) \in V_{4} \) we consider backward iterations of \( f_{a}  \).
Note that \( (x_{-1}, y_{-1})= (y/b, x-1+ay^{2}/b^{2})
+\tilde \varphi  (x,y,a)\).
Then \( x_{-1}\geq 2/b -\eta \geq -2 \) and \( y_{1}\geq
-2+4a/b^{2}-\eta \geq 2 \). Thus \( f^{-1}(x,y)\in V_{4} \) and \(
y_{1}\geq y/b \). Therefore \( y_{-n}\geq y/b^{n} \) and so \(
|y_{-n}|\to\infty \). 

For \( (x,y) \in V_{5} \) we have \( y\geq |x| \geq 2 \). Thus \(
x_{1}\geq y/b-\eta  \geq 2 \) and \( y_{1}\geq y^{2}/b^{2}\geq 2 \) .
So we have that \( f^{-1}(x,y)\in V_{4} \), and we argue as above.

For \( (x,y) \in V_{6} \), we have \( x_{-1}=y/b+\tilde \varphi
_{2}(x,y,a) \geq 2 \) and \(
y_{-1}\geq 2 \). Therefore, \( f_{a}(x,y) \in V_{4} \) and again, we
argue as above.

    \end{proof}

    \begin{lemma}
\( \mathcal D  \subset \hat R\). 
\end{lemma}

\begin{proof}
  The arguments used above have implications 
for the locations of the stable and unstable manifolds of the fixed 
points. Indeed the stable  manifolds of the fixed points 
cannot intersect \( V_1\cup V_2 \cup V_3 \) since all points in this 
region tend to infinity in forward time, whereas, by definition, 
points in the stable manifolds tend to the fixed points under forward 
iteration. Similarly the unstable 
manifolds of the fixed points cannot intersect \( V_{4}\cup V_{5}  \)
since all points in this region tend to infinity in backwards time. 
By definition \( \mathcal D \) is bounded by arcs of stable and unstable
 manifolds of the fixed point as in the Figure
 and therefore \( \mathcal D \subset \mathbb R^{2}\setminus
(V_{1}\cup \ldots \cup V_{5}) = (-2, 2)\times (-4, 2) \).
    \end{proof}
\begin{figure}
       \includegraphics[width=9cm]{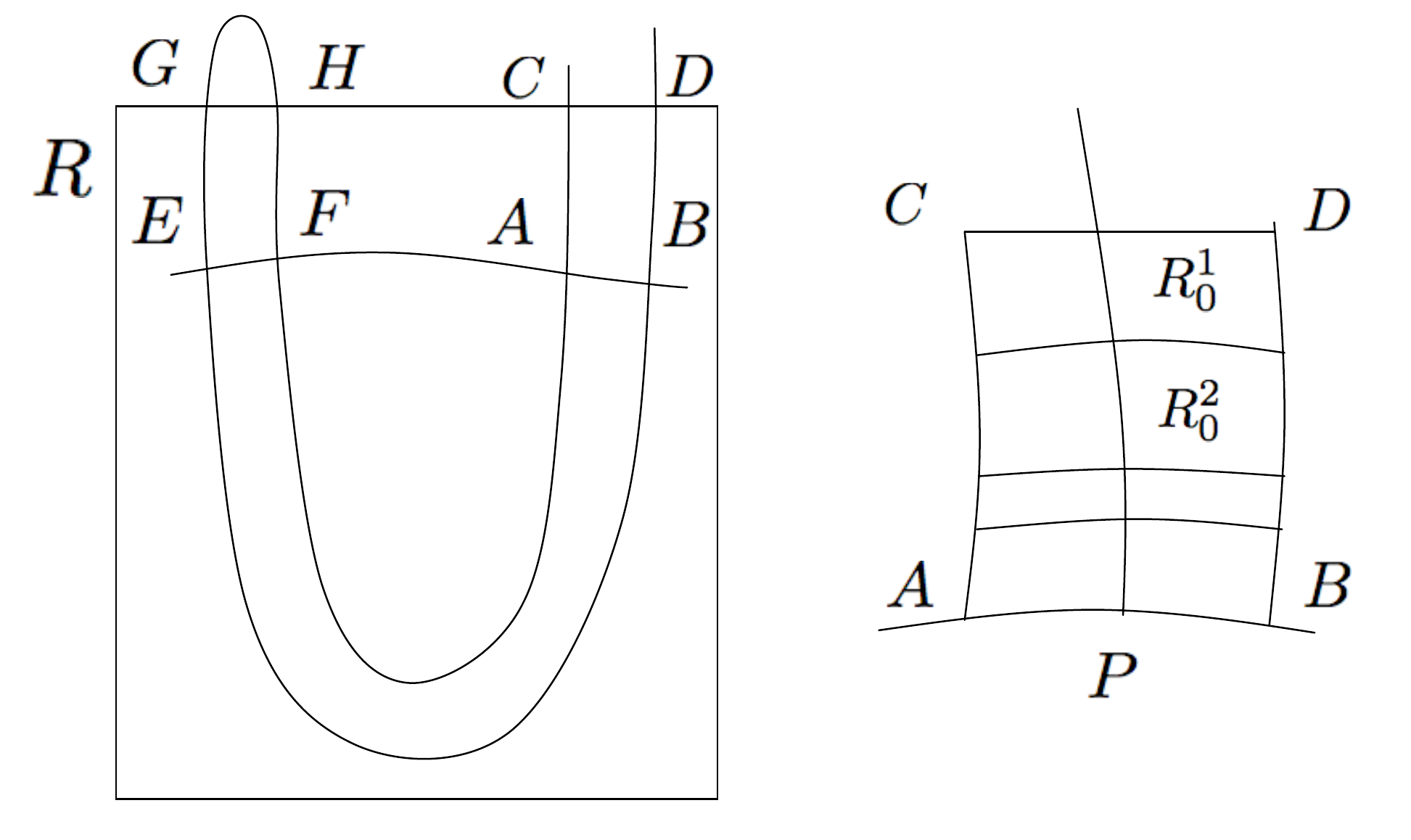}
   \caption{Regions \(R_{0}^{i}\)}\label{regionsR}  
   \end{figure}

 \begin{lemma}
\( \Omega\subset \mathcal D. \)
\end{lemma}
   
\begin{proof}
   To show that \( \Omega\subset \mathcal D \) 
   we refine the strategy used in the proof of the previous 
    lemma, and show that the orbits of all points outside \( \mathcal 
    D \) are unbounded in either backward or forward time. Since we
    have already shown that \( \Omega\subset R\), we need 
    to consider only points in the region \( R \setminus \mathcal D \). 
    
    \subsubsection*{Subdividing}
  We write 
  \[
  R \setminus \mathcal D = R_{0}\cup R_{1}\cup R_{2}\cup
  R_{3}
  \]
  where the regions 
  \(  R_{0}, R_{1}, R_{2}, R_{3}\)   are defined as follows. 
  Consider the points \(A\) and \(B\) of intersection of
     \(W^{u}(P)\) and \(W^{s}(q)\) and \(C\) ,\(D\) of intersection of 
     \(W^{s}(q)\) and \(y=4b\) as in the Figure \ref{regionsR}. 
 We let   \(R_{0}\) denote the closed region
     bounded by the arcs of manifold \(AC\), \(AB\) and \(BD\), and the
     segment \(CD\). 
  We let \( R_{1} \) denote   
     the
     region bounded by the arc of manifolds  \( HF\), \(FA\) and \(AC\), and 
     the segment \(HC\).
  We let \( R_{2} \) denote the region  bounded by the arcs of
  \(W^{u}(p)\) and \(W^{s}(q)\) between the points \(E\) and \(F\), as
  in Figure \ref{regionsR}. 
  We let \( R_{3} = R \setminus (\mathcal D \cup R_{0}\cup R_{1}\cup
  R_{2} )\).  
  We also define
  \[ \tilde R_{3}\subset R_{3}
  \]
  as  the region satisfying \(-2b-\eta <y<2b+\eta\) at the
left side of the arc of \(W^{s}(q)\) between the points \(I\) and
\(J\), of intersection of that manifold with the lines \(y=-2b-\eta\) 
and \(y=2b+\eta\), as in the Figure \ref{RegionsRNew} (b).

     \subsubsection*{Points of  \protect\( R_{0} \protect\)
     escape in backward time}
   Since \(b\) is small, we have that all the points 
    \((x,y) \in R_{0}\) satisfy \(x>0.2\).  
Notice that, for the unperturbed H\'enon map
    \(h_{a,b}(x,y)=(1-ax^{2}+y, bx)\), we have that
     any piece of curve  \( \gamma \) with slope 
    less than \( 1/10 \) contained in the region where \(|x|>0.2\)
    is mapped to another curve with slope less than \( 1/10 \). 
    Indeed, 
    letting \( (v_{1}, v_{2}) \) denote a tangent vector to \( \gamma \) 
    with \( |v_{2}|/|v_{1}| < 1/10 \), we have \( (v_{1}', v_{2}') = 
    Dh_{a,b}(v_{1}, v_{2}) = (-2axv_{1}+v_{2}, bv_{1}) \) whose slope is 
    \( |v_{2}'|/|v_{1}'| = |b/(-2ax+(v_{2}/v_{1}))|<1/10 \), provided 
    \(b\) is small and \(a\) is close to 2. For future reference,
    notice that, if \(|x|>0.5\), we also have that the norm of
    \((v_{1},v_{2})\) is uniformly expanded. So, since \(f_{a}\) is close to  
    \(h_{a,b}\) in the \(C^{2}\) topology, we can assume that \(f_{a}\) 
    also has this property in \(R_{0}\).
    
    Now call \(\alpha_{n}\) the successive images of the segment \(CD\)
    intersected to \(R_{0}\). Since they cannot intersect each other, 
    and \(CD\) has a point of the stable manifold of \(p\), the curves
    \(\alpha_{n}\) determine a system of ``fundamental domains'' in
    \(R_{0}\): they cross \(R_{0}\) from one stable boundary to the
    other, and they converge to the arc of unstable manifold \(AB\). 
    Call \(R_{0}^{i}\) the region of \(R_{0}\) between 
    \(\alpha_{i-1}\) and \(\alpha _{i}\), \(\alpha _{0}=CD\), and notice that 
    \(f^{-1}(R_{0}^{i})\subset R_{0}^{i-1}\). We also have that
    \(f^{-1}(R_{0})\) falls outside \(R\). That implies that
    \( R_{0} \setminus AB\)  does not intersect \(\Omega \), and any point
    which has 
    an iterate in \(R_{0}\setminus AB\) is not in \(\Omega\).

    \subsubsection*{Points of  \protect\( R_{1} \protect\)
    map to \protect \( R_{3} \protect\)}
We show that \( f(R_{1}) \cap R \subset R_{3} \). 
Indeed, the unstable eigenvalue of \( p \) is positive and therefore 
\( f(R_{1}) \) must remain on the same side of \( W^{s}(q) \) as \(
R_{1} \). Also, since \( f(R) \subset \mathbb R \times [-2b-\eta, 
2b + \eta]  \) we have that \( f(R_{1}) \) does not intersect any of \(
\mathcal D, R_{0}, R_{1}, R_{2} \).  

\begin{figure}
    \includegraphics[width=9cm]{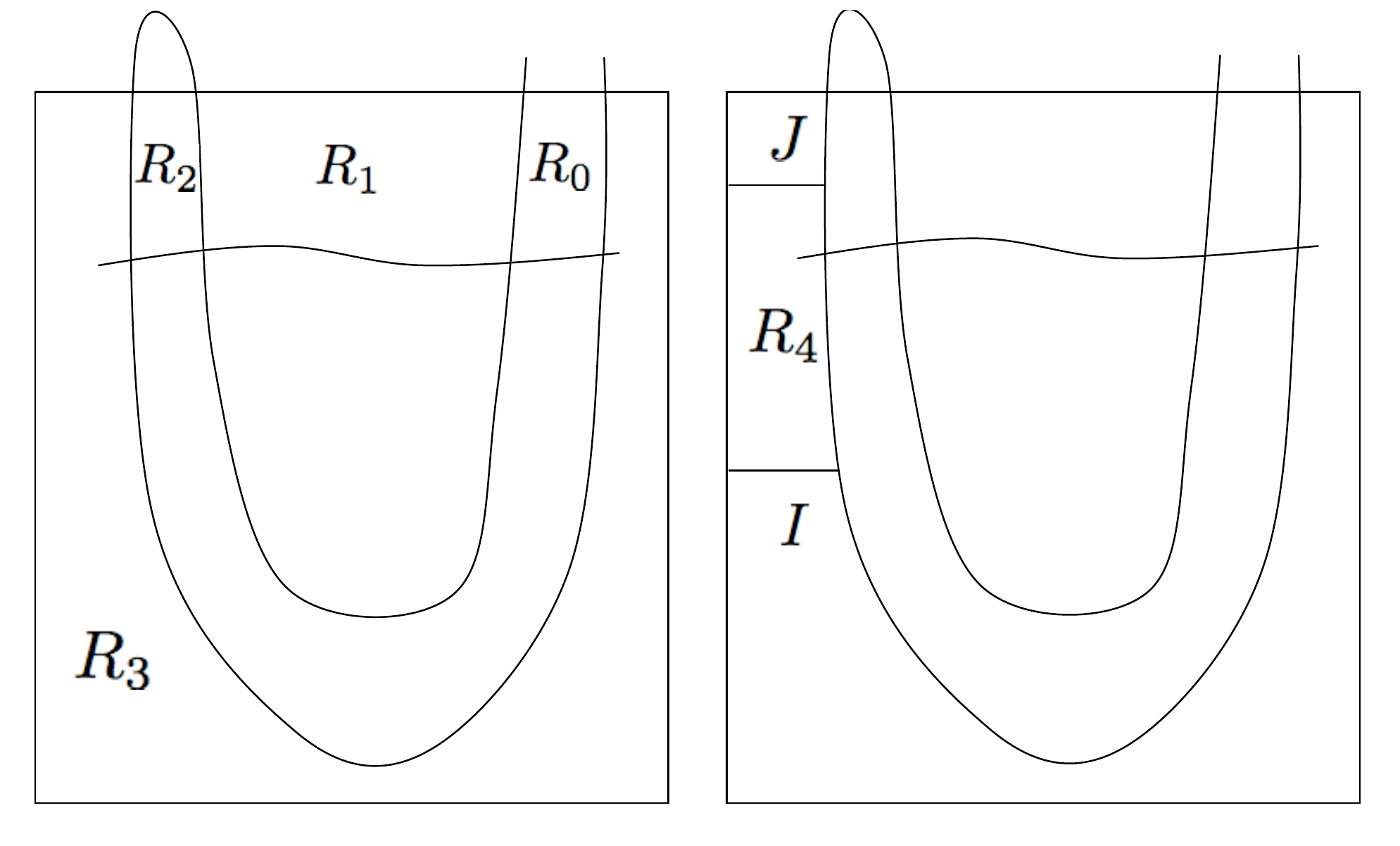}
    \caption{Regions \( R_{1} \) to \( R_{4} \)}\label{RegionsRNew}
 \end{figure} 

 \subsubsection*{Points of  \protect\( R_{3} \protect\)
 map to \protect \( \tilde R_{3} \protect\)}
 We
 now show that \( f(R_{3}) \subset \tilde R_{3}\).  
Again, we use the fact that \( f(R) \subset \mathbb R \times [-2b-\eta, 
2b + \eta]  \). Then,  since one of the components of the
boundary of \(R_{3}\) is an arc of stable manifold of \(q\)
containing the fixed point \(q\), and the unstable eigenvalue of \(q\)
is positive, we conclude that the image of
\(R_{3} \) is contained \( \tilde R_{3} \).

\subsubsection*{Points of  \protect\( \tilde R_{3} \protect\)
escape in forward time}

We can assume, if
\(b\) is small, that all the points \((x,y)\) in \(\tilde R_{3}\) satisfy
\(x<-0.5\) (notice that, for \(b=0\), we have \(q=(-1,0)\). 
Take \(t\) a point in \(\tilde R_{3}\setminus W^{s}(q)\), and 
connect \(t\) to the
boundary of \(\tilde R_{3}\) through a horizontal line inside \(\tilde
R_{3}\), determining a point \(t' \in W^{s}(q)\).
Again, by the proximity of \(f\) and \(h_{a,b}\), and the fact that 
vectors with slope smaller than \(1/10\) in \(\tilde R_{3} \cap R\)  are sent by
\(Dh_{a,b}\) in vectors with slope smaller than \(1/10\), and uniformly 
expanded, we have that the horizontal distance between \(f(t)\) 
and \(f(t')\) is uniformly expanded. Applying \(f\) repeatedly, as
long as the image is inside \(\tilde R_{3}\cap R\), we have that the
horizontal distance between the successive images of \(t\) and
\(W^{s}(q)\) increases exponentially. Then the forward images of \(t\) 
leave \(R\) for some positive time.

\subsubsection*{Points of  \protect\( R_{2} \protect\)
map to \protect\( R_{0} \protect\) in backward time}

Notice that \( f^{-1}(R_{2})\cap R \subset R_{0} \) since 
all the other 
regions in \(R\) outside \(\mathcal D\) are mapped forward to the region
\(\tilde R_{3}\), and so do not contain points of the backward image of
\(R_{2}\). Moreover, the unstable boundary of \(R_{2}\) belongs to
\(W^{u}(p)\) and approaches \(p\) as we apply \(f^{-1}\), and the
stable boundary cannot cross \(W^{s}(q)\), then \(f^{-1}(R_{2})\)
does 
not intersect \(\mathcal D\). 
Since \(f^{-1}(R_{2})\cap R\subset R_{0}\), the points
in there that are not 
in \(W^{u}(p)\) leave \(R\) for backward iterations. 
\end{proof}

\begin{figure}
    \includegraphics[width=9cm]{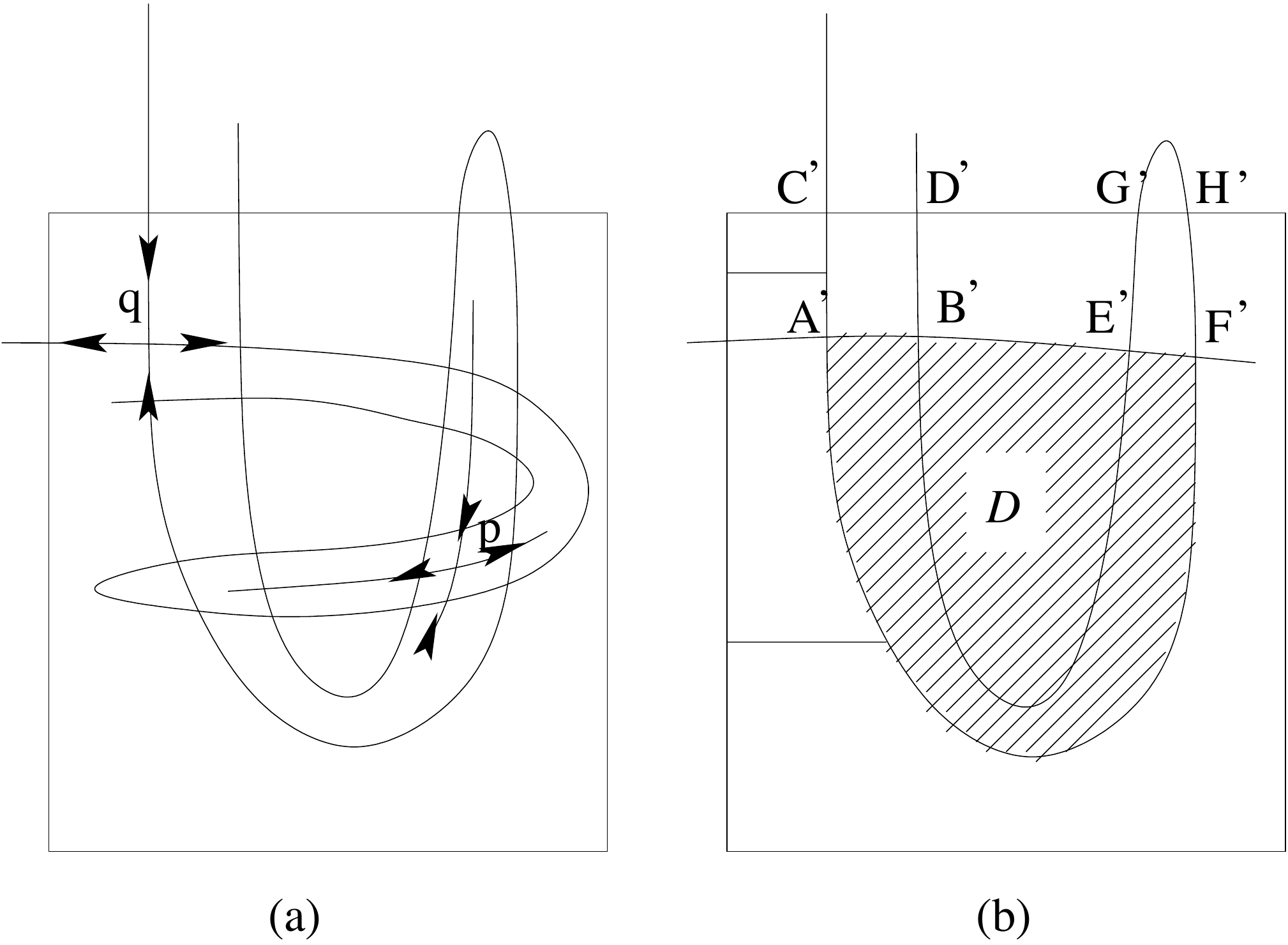}
    \caption{Invariant manifolds and the region \( \mathcal D \) for 
    \( b< 0 \)}\label{negativeb}
 \end{figure}

\begin{lemma}
\( \Omega \subset \overline{W^{u}(p)} \). 
\end{lemma}
\begin{proof}
Notice first of all that by the \( \lambda \)-lemma we have \( 
q\in\overline{W^{u}(p)} \). Now suppose by contradiction that there 
exists \( z=(x,y) \in \Omega \) with \( z\notin \overline{W^{u}(p)} \). 
Then there exists \( \varepsilon \) and an \( \varepsilon \) 
neighbourhood \( B_{\varepsilon}(z) \) of \( z \) with \( 
B_{\varepsilon}(z) \cap \overline{W^{u}(p)} = \emptyset \). 
Since \( \Omega \) is \( f \)-invariant we have \( 
f^{-n}(z)\in\Omega(f)\subset\mathcal D \) for all \( n\in\mathbb N \) 
and therefore \( z\in f^{n}(\mathcal D) \) for all \( n\in\mathbb N 
\). Notice that the boundary \( \partial f^{n}(\mathcal D) \subset 
W^{u}(p) \cup  f^{n}(EB^{s}) \) (where \( EB^{s} \) denotes the piece of \( 
W^{s}(q) \)  between \( E \) and \( B \), as in Figure \ref{regionsR}. 
It is enough therefore to show that, for large \(n\), the \( \partial 
f^{n}(\mathcal D) \) is \( \varepsilon \)-dense in \( f^{n}(\mathcal D) 
\) as this will imply that \( B_{\varepsilon}(z) \cap 
\overline{W^{u}(p)} \neq \emptyset  \), contradicting the assumptions. 
This follows easily from the fact that \( f \) is (strongly) area 
contracting, and thus the area of \( f^{n}(\mathcal D)  \) tends to 
zero as \( n\to\infty \). In particular we must have that \( 
B_{\varepsilon}(z)\cap \partial f^{n}(\mathcal D) \neq \emptyset\) for 
all \( n\geq N \) sufficiently large. 
Moreover, the length of the part 
of the boundary which belongs to \( W^{s}(q) \) also tends to zero. 
Thus most of the boundary belongs to \( W^{u}(p) \) and thus we must 
have \( 
B_{\varepsilon}(z)\cap W^{u}(p) \neq \emptyset\) for all \( n \) 
sufficiently large.
\end{proof}

\subsubsection{Completion of the proof}
\begin{proof}[Proof of Proposition \ref{nonwan}]
 Combining the results of the Lemmas stated above we have that 
 \( \Omega\subset  \overline{W^{u}(p)} \cap \mathcal D\).  The
 statement in the Proposition now follows immediately by observing that
 \( \Omega \subset \mathcal D \) implies \( \Omega \subset f(\mathcal 
 D) \) and that 
 \( f(\mathcal D) \subset [-2,2]\times [-4b, 4b] \) directly from the 
 definition of \( f \) if \( \eta \) is sufficiently small.

Finally,  in the case \(b<0\), we consider
the stable and unstable manifolds of
\(q\) crossing as in Figure \ref{negativeb} (the rectangle \(R\) is
exactly the same), determining the region \(
\mathcal D\) in this case. The proof is entirely analogous considering
the points \( A', B' \), etc., corresponding to the points \( A, B, \)
etc., above. 
    \end{proof}

\section{Hyperbolic coordinates and critical points}
\label{sectionhyp}

The key idea of our whole strategy is the notion of \emph{dynamically
defined critical point} which relies on the fundamental notion of
\emph{hyperbolic coordinates}. In this section we introduce these
notions and develop
the main technical ideas which we will use. 
In Section \ref{prelim} we clarify the relations between various
constants used in the argument and introduce some preliminary
geometric constructions. In Section \ref{hypcoord} we discuss the
definition and basic theory of hyperbolic coordinates. In Section 
\ref{curvature} we introduce the idea of \emph{admissible curves} and 
prove certain estimates concerning the images of admissible curves.
Finally, in Section \ref{hypcrit} we introduce the notion of
dynamically defined critical point and prove that such critical points
always exist in images of certain admissible curves.

\subsection{Preliminary geometric definitions and fixing the constants}
 \label{prelim}

\subsubsection{Fixing the constants}\label{constants}
We now explain  the required relations between the different constants
used in the proof, and the order in which these constants are chosen. 
All constants are positive.
First of all we fix two constants 
\[ 
\delta = 1/10 \quad\text{ and } \alpha = 1/2.
\]
These will be introduced in Sections \ref{fixptnhbd} and
\ref{curvature} below. Even though we specify the actual numerical value
of these constants we shall continue to use the constants in the
argument because they have some specific geometric meaning and it is
useful to keep track of their occurrence throughout the paper. 
We then  fix a constant 
\( k_{0} \) large enough so that
\begin{equation}\label{k0}
\frac{\sqrt{\delta}}{2\sqrt{3}} 
\left(\sqrt{3/\sqrt 5}\right)^{k_{0}-1}>1 
\end{equation}
In Section \ref{critnhbd} we fix a constant \( \varepsilon \) which
will then remain unchanged. 
Finally, at some finite number of places in the argument, 
we will require \( a \) to be sufficiently close to 2 and \( |b| \)
and \( \eta \) to be sufficiently small. 

We remark that we can suppose  that \( a \) 
is close to 2 without compromising the fact that hyperbolicity holds
for all larger values of \( a \). Indeed, once we fix a neighbourhood 
of \( 2 \) in the \( a \) parameter space, we can always guarantee
uniform hyperbolicity for values of \( a> 2 \) outside this
neighbourhood by taking \( |b| \) and \( \eta \) sufficiently small,
(depending on this neighbourhood).

\subsubsection{The fixed point neighbourhoods}
\label{fixptnhbd}
\begin{figure}
    \includegraphics[width=9cm]{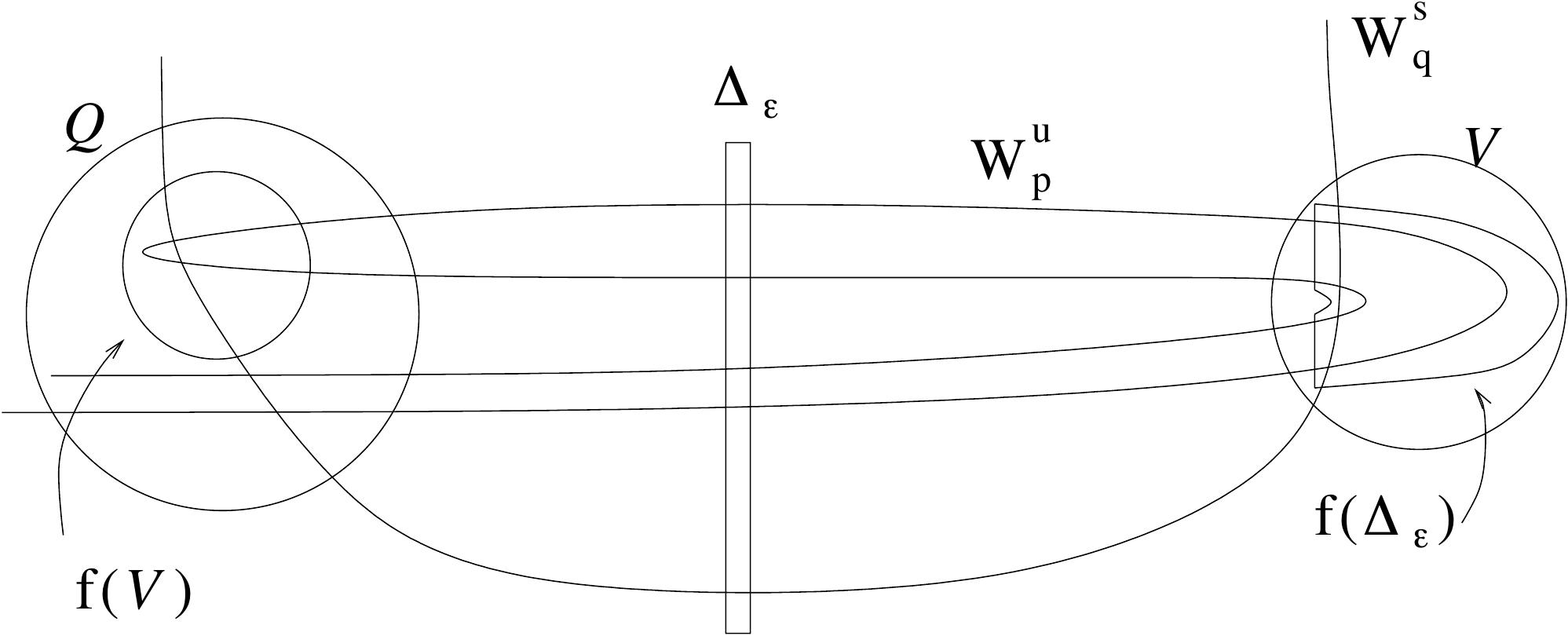}
  \caption{The neighbourhoods \(\mathcal Q\) and \(\mathcal V\)}
  \label{Q}
  \end{figure}

Recall first of all that the map \( f_{*}=h_{2,0} \) 
has two fixed points \( p_{*} \) 
and \( q_{*}=(-1,0)  \) with  \( f_{*}(1,0) =   q_{*} \). 
For \( \delta=1/10 \) we let  
\( \mathcal Q = \mathcal Q_{0}:= B_{\delta}(q_{*}) \) 
be the open ball of radius \( \delta \) 
centred at \( q_{*} \) and \( \mathcal V=\mathcal V_{0} \) be the
component of \( f^{-1}_{*}(\mathcal Q) \) not intersecting \( \mathcal Q
\), see Figure  \ref{Q}. 
The expanding eigenvalue at the point \( q_{*} \) is equal to 4 and so
we can suppose that \(|a-2|,  |b|, \eta, \) are all 
small enough so that \( \|Df_{z}\| > 3\) for all \( z\in \mathcal Q \).
Then, for \( n\geq 0 \), let 
\[ 
\mathcal Q_{n}(f)=\bigcap_{i=0}^{n}f^{-i}(\mathcal Q_{0}) \quad \text{ and }
\quad \mathcal V_{n}(f)=f^{-1} (\mathcal Q_{n}(f)) \cap \mathcal V.
\]
Notice that \( \mathcal V_{n} \) 
is just the component of \( f^{-1} (\mathcal Q_{n}(f)) 
\) containing \( (1,0) \).  
Since \( \mathcal Q_{n} \) is a 
neighbourhood of \( q \) for every \( n \), \( \mathcal V_{n}\setminus
f^{-1}(W^{s}_{\delta}(q))\), 
where \( W^{s}_{\delta}(q) \) denotes the connected component
containing \( q \) of \( W^{s}(q) \cap \mathcal Q_{0} \),
has two components: 
we let 
\[ 
\mathcal V_{n}^{-}=\mathcal \mathcal V_{n}\cap \mathcal D 
\quad \text{ and }\quad 
\mathcal V_{n}^{+} = \mathcal V_{n} \setminus \mathcal V_{n}^{-}.
\]
Notice that a piece of  \( W^{s}(q) \) forms the boundary between 
\( \mathcal V_{n}^{+} \) and \( \mathcal V^{-}_{n} \). 
We mention for
future reference a 
 simple estimate which we shall use below. 
\begin{lemma}\label{distance1}
\( 
    d(z, f^{-1}(W^{s}_{\delta}(q)) ) \geq \delta/5^{k} 
\) for all  \( z\in \mathcal
    V_{k}\setminus  \mathcal V_{k+1}.\)
 \end{lemma}
 \begin{proof}
\( z \in \mathcal V_{k}\setminus\mathcal V_{k+1} \) implies, by
definition,  \( d(z_{k+1}, q ) \geq \delta \). 
For points \( z \) close to \(
f^{-1}(W^{s}_{\delta}(q) ) \) 
this also means \( d(z_{k+1}, W^{s}_{\delta}(q))
\geq \delta\) since such points come very close to the fixed point \( 
q \) and escape the \( \delta \) neighbourhood of \( q \) along the
direction of \( W^{u}(q) \). Thus, using the fact
that the norm of the derivative \( Df \) in \( \mathcal D \) is
uniformly bounded above by 5 we obtain the
result. 
\end{proof}

\subsubsection{The critical neighbourhood}
\label{critnhbd}
For \( \varepsilon > 0 \) we define a 
\emph{critical neighbourhood} 
\[
\Delta_{\varepsilon} = (-\varepsilon, \varepsilon) 
\times (-4b, 4b).  
\]
Notice that we can take \( \varepsilon \) sufficiently small so that 
\( q_{f}\in f(\mathcal V) \) and 
\[ 
f(\Delta_{\varepsilon}) \subset \mathcal V_{k_{0}}.
\]
From now on we consider \( \varepsilon \) fixed. 
We also let 
\begin{equation}\label{delta}
\Delta=\Delta_{a}=\{x\in\Delta_{\varepsilon}: f (x) \notin\mathcal 
D\}.
\end{equation}
For \( a \) sufficiently close to \( 2 \) and \( |b| \) and \( \eta \)
sufficiently small we have uniform hyperbolicity outside \(
\Delta_{\varepsilon} \).  We state this fact more formally in the
following
\begin{lemma} \label{LemmaUnifHyp}
For every 
\( \hat\lambda \in (0, \log 2) \) and 
\( |a-2|, |b|, \eta >0 \) sufficiently small,
there exists a constant \( C_{\varepsilon}>0 \) 
such that for all \(  k\geq 1 \) and points 
\( z \) with   \( z, f(z), \ldots, 
 f^{k-1}(z) \notin \Delta_{\varepsilon} \), and vector \( v  \) with 
 slope \(  <  \alpha \) we have 
\begin{equation}\label{slope}
 \text{slope }  Df^{k}_{z}(v)   < \alpha,  
 \end{equation}
 \begin{equation}\label{UE1} 
 \|Df^{k}_{z}(v)\|\geq C_{\varepsilon} e^{\hat\lambda  k}\|v\|.
 \end{equation}
  If,  \emph{moreover}, \( f^{k}(z)\in\Delta \) 
 then we have 
   \begin{equation}\label{UE2}
 \|Df^{k}_{z}(v)\|\geq e^{\hat\lambda k}\|v\|.
 \end{equation}
 \end{lemma}
 \begin{proof}
This is a standard result (see for example \cite{BenCar91} or
\cite{MorVia93}) and so we omit the details. We just mention that it
follows from the fact that the limiting one-dimensional map 
\( h_{2,0} \) satisfies uniform expansivity estimates outside an
arbitrary critical neighbourhood \( \Delta_{\varepsilon} \) (with
constant \( \hat\lambda \) arbitrarily close to \( \log 2 \) but
constant \( C_{\varepsilon} \) depending on \( \varepsilon \) and arbitrarily small
for  \(  \varepsilon \) small), see e.g. \cite{MelStr93}. 

Considering 
this one-dimensional map as embedded in the space of two-dimensional
maps and 
using the fact that uniform hyperbolicity is an open condition
we obtain the statement in the Lemma for \( |b|, \eta \neq 0 \) 
sufficiently small. 
\end{proof}

\subsection{Hyperbolic coordinates}  
\label{hypcoord}       

The notion of Hyperbolic Coordinates is inspired by 
some constructions in \cites{BenCar91, MorVia93},
developed in  \cite{LuzVia03} and formalized in 
\cite{HolLuz06} as an alternative framework with 
which to approach the classical theory of invariant manifolds.  
Here we review the basic definitions and theory to the extent to which
they will be required for our purposes.

\subsubsection{Hyperbolicity of compositions of linear maps}
We recall the notion of hyperbolic coordinates and give the basic
definitions and properties in the general context of \( C^{2} \)
diffeomorphisms of a Riemannian surface \( M \). 
For \( z\in M,  k\geq 1 \) we let 
\[ 
 F_{k}(z) = \|Df^{k}_{z}\|\quad\text{ and } \quad E_{k}(z) =
 \|(Df^{k}_{z})^{-1}\|^{-1} 
\]
denote the maximum expansion and the maximum contraction respectively 
of \(  Df^{k}_{z}\). Then we think of the quantity
\[ 
H_{k}(z) =
 E_{k}(z) / F_{k}(z)
 \]
as the \emph{hyperbolicity} of \(  Df^{k}_{z}\). 
Notice that \( H_{k}\leq 1 \) always. 
The  condition \( H_{k} = E_{k}/F_{k} < 1 \)
implies that the linear map \( Df^{k} \) maps 
the unit circle \(\mathcal S \subset T_{z}M \) 
to an  ellipse \( \mathcal S_{k} =
Df_{z}^{k}(\mathcal S) \subset T_{f^{k}(z)}M \)
with well defined major and minor axes. 
The unit vectors \( e^{(k)}, f^{(k)} \) 
which are mapped to the minor and major
axis respectively of the ellipse, and are thus the \emph{most 
contracted} and \emph{most expanded} vectors respectively,  are given
analytically as solutions to
 the differential equation \( d\|Df_{z}^{k}(\cos\theta, 
 \sin\theta)\|/d\theta = 0 
 \) which can be solved to give the explicit formula
\begin{equation}\label{contdir}
\tan 2\theta  =
\frac{2 [(\pfi x1^{k})(\pfi y1^{k}) +(\pfi x2^{k})(\pfi y2^{k})]}
{(\pfi x1^{k})^2+(\pfi x2^{k})^2 - (\pfi y1^{k})^2 -(\pfi y2^{k})^2}.
\end{equation}
Here \( f=(f_{1}, f_{2}) \) are the two coordinate functions of \( f \). 
Notice that \( e^{(k)} \) and \( f^{(k)} \) are always \emph{orthogonal}
 and \emph{do not} in general correspond to the 
stable and unstable eigenspaces of \( Df^{k} \). 

\subsubsection{Hyperbolic coordinates and stable and unstable
foliations}

We define the 
\emph{hyperbolic coordinates of order} \( k \) at the point \( z \) 
as the orthogonal coordinates \( \mathcal H_{k}(z) \) 
given at \( z \) 
by the most contracted and most expanded 
directions for \( Df^{k}_{z} \).  
If \( f \) is \( C^{2} \) and \( H_{k}(z) < 1 \) then hyperbolic
coordinates are defined in some neighbourhood of \( z \) and define
two orthogonal \( C^{1} \) vector fields. In particular they are
locally integrable and thus give rise to 
two orthogonal  foliations.  
We let \( \mathcal E^{(k)} \) denote the 
\emph{stable foliation of order} \( k \) 
formed by the integral curves of the vector field \( \{e^{(k)}\} \) 
and \( \mathcal F^{(k)} \) denote the \emph{unstable foliation of 
order} \( k \) formed by the integral curves of the vector field \( 
\{f^{(k)}\}\).

\subsubsection{Hyperbolic coordinates for H\'enon-like maps}
\label{hyphen}

A crucial property of hyperbolic coordinates and finite order
stable and unstable foliations is that, under very mild assumptions, 
they \emph{converge} in quite a strong sense as \( k\to \infty \). 
We formulate a version of this property here in our
specific context.

\begin{proposition}\label{c2close}
For every \( k\geq 1 \), 
hyperbolic coordinates \( \mathcal H_{k} \) and stable and unstable
foliations \( \mathcal E^{(k)} \) and \( \mathcal F^{(k)} \) are
defined in \( \mathcal V^{+}\cup \mathcal V^{-}_{k} 
\) 
Moreover 
\begin{enumerate}
    \item the angle between each stable direction \( e^{(k)} \) and 
    the slope of \( f^{-1}(W^{s}_{\delta}(q)) (\approx 2) \).
\item the curvature of each stable leaf ,
\end{enumerate} 
are both \( \lesssim b \). 
Also, the \( C^{2} \) distance between 
leafs of \( \mathcal E^{(k)} \)
and leaves \( \mathcal E^{(k+1)} \) is 
\( \lesssim bk \). 
\end{proposition}

\begin{proof}
    
Analogous
convergence results are formulated and proved in great generality in
\cite{HolLuz06}  under weak (subexponential) growth of the
derivative. 
    Here we shall need only some very particular cases of
    these estimates and therefore we first describe the specific setting
    in which they will be applied here.
 The main ingredient for the proof is that fact that    
by our choice of \( \delta \) and assuming that \( |a-2| \), \(
 |b| \) and \( \eta \) are small enough we have that
 \( \|Df (z) - Df_{*}(q_{*})\| \) is small for all \( \ z\in \mathcal Q \)
 and thus in particular 
 \begin{equation}\label{contest}
  E_{k}(z_{0}) \leq b^{k} \quad\text{and} \quad 
  F_{k}(z_{0}) \geq 3^{k} \quad \forall \ z_{0}\in \mathcal V_{k}. 
  \end{equation} 
It follows  immediately that hyperbolic coordinates, and their
associated foliations, of order \( k \) 
are well defined in \( \mathcal V_{k} \). 
Points in \(  \mathcal V_{k}^{-} \) are then re-injected into 
\( \mathcal D \setminus \mathcal Q \) and these hyperbolicity
estimates can no longer be guaranteed, a priori, for all time. Points 
in \( \mathcal V_{k}^{+} \) however are outside \( \mathcal D \) 
and therefore, by the arguments of Section
\ref{proofnonwan}, eventually escape towards infinity. In particular
the required hyperbolicity conditions can be guaranteed to hold for
all positive iterates. This implies that hyperbolic coordinates of
order \( k \) are well defined in \( \mathcal V^{+}\cup\mathcal V^-_{k} \)
as in the statement of the proposition.

The statements about the direction of the stable directions, the
curvature of the leaves and the \( C^{2} \) distance between stable
leaves of different orders, all follow directly from 
\cite{HolLuz06}*{Main Theorem}. These calculations are purely technical
and do not add to our geometrical understanding of this situation, we 
therefore omit the details and refer the reader to that paper.

\end{proof}

\subsection{Admissible curves}
\label{curvature}

Recall that the curvature \( \kappa(s) \) of a parametrized curve 
\( \gamma(s) = (x(s), y(s)) \) is given by 
\[ 
\kappa(s) = \frac{|\dot x \ddot y-\dot y \ddot x|}{\|(\dot x, \dot
y)\|^{3}} =
\frac{|\dot\gamma \times \ddot\gamma |}{|\dot\gamma|^{3}}.
\]
The equivalence between the two formulas is given by the formula
 \((v_{1},v_{2})\times
(w_{1},w_{2})=v_{1}w_{2}-v_{2}w_{1}\). 
\begin{definition}
   For  \( \alpha > 0 \), we say that  a 
   \( C^{2} \) curve \( \gamma = \gamma(s) = (x(s), y(s)) \) is 
\emph{admissible}  if 
   \(    |\dot y(s)|/|\dot x(s)| < \alpha \) and     
   \( |\kappa(s)| < \alpha \) for all \( s \).  
    \end{definition}
We remark that both the curvature and the slope of tangent vectors of 
a curve are independent of the parametrization, and thus so is the
definition of admissibility. 
We shall want to compare the curvature at a point of a curve and at
the corresponding point of its image. So, we suppose \(
\gamma_{i-1}(s) \) is a parametrized \( C^{2} \) curve and 
\( \gamma_{i}(s) =f(\gamma_{i-1}(s)) \). 
For simplicity we shall often omit the parameter \( s \) and
simply write \( Df \) to denote the derivative at the point \(
\gamma_{i-1}(s) \). 

\begin{proposition}\label{smallcurv}
    Let 
    \( \{\gamma_{i}\}_{i=0}^{n}\) be a sequence of  \( C^{2} \)
   curves with  \( \gamma_{i}=f^{i}(\gamma_{0}) \). 
Suppose that for some \( s \), \( n \) is a ``hyperbolic time'' in the sense that 
\[
  \|\dot\gamma_{n}(s)\|\geq C e^{\lambda 
  j}\|\dot\gamma_{n-j}(s)\| 
\] 
for all \( j=0, .., n-1 \).
Then for \( |b|, \eta \)
sufficiently small,  \(
  \kappa_{0}(s) < \alpha \) implies \(\kappa_{n}(s) < \kappa_{0}(s) < 
  \alpha. \)
  \end{proposition}

    \begin{corollary}
	\label{admisstoadmiss}
 If \( \gamma\subset \mathcal 
	D\setminus\Delta_{\varepsilon} \) is admissible, then \( 
	f(\gamma) \) is also admissible
      \end{corollary} 

      \begin{proof}
This follows from Proposition  \ref{smallcurv} and 
the  hyperbolicity outside \( \Delta_{\varepsilon} \).
Condition \eqref{slope} implies that the slope of each tangent
vector to \( f(\gamma) \) is \( <\alpha \) and condition \eqref{UE1}
together with Lemma \ref{smallcurv} gives the curvature \( <\alpha \).
   \end{proof}

To prove Proposition \ref{smallcurv} we first 
prove a general curvature
estimate. We fix some bounded 
neighbourhood \( \hat R \) of \( R \) and, as above,   
we suppose 
\( \{\gamma_{i}\}_{i=0}^{n}\) is a sequence of  \( C^{2} \)
   (not necessarily admissible) curves with  \( \gamma_{i}=f^{i}(\gamma_{0}) \), 
all contained in \( \hat R\).

\begin{lemma}\label{curv}
There exists  $K > 0 $ independent of \( a, b, \eta \) 
such that for all \( i=1,\ldots, n \) we have 
\[
\kappa_i(s) \le K(b+\eta) 
\frac{|\dot{\gamma}_{i-1}(s)|^3}{|\dot{\gamma}_{i} (s)|^3} 
\kappa_{i-1}(s)+ 
K(b+\eta) \frac{|\dot{\gamma}_{i-1} (s)|^3}{|\dot{\gamma}_{i} (s)|^3} 
\]
\end{lemma}
\begin{proof}
    We use the formula 
    \( \kappa = |\dot\gamma \times \ddot\gamma |/|\dot\gamma|^{3} \)
    for the curvature. 
We have 
\[ 
\dot\gamma_{i}= (Df) \dot \gamma_{i-1}   = 
\begin{pmatrix} f_{1,x} & f_{1,y} \\  f_{2,x} & f_{2,y}
\end{pmatrix}\dot \gamma_{i-1} 
= \begin{pmatrix} -2ax_{i-1} + \varphi_{1,x}
    &1+\varphi_{1,y}\\
    b+ \varphi_{2,x} & \varphi_{2,y}
    \end{pmatrix} \dot \gamma_{i-1} 
\]
and 
\[ 
\ddot\gamma_{i} = \begin{pmatrix} 
\nabla f_{1,x} \cdot \dot\gamma_{i-1} & \nabla f_{1,y} \cdot 
\dot\gamma_{i-1} 
\\ \nabla f_{2,x} \cdot \dot\gamma_{i-1} & \nabla
f_{2,y}\cdot \dot\gamma_{i-1} 
\end{pmatrix}
\dot\gamma_{i-1} +
(Df)\ddot\gamma_{i-1}.
\]
Therefore \( \dot\gamma_{i}\times\ddot\gamma_{i}  \) is given by 
\begin{equation}\label{curv1}
(Df)\dot\gamma_{i-1} \times 
\begin{pmatrix} 
\nabla f_{1,x} \cdot \dot\gamma_{i-1} & \nabla f_{1,y} \cdot 
\dot\gamma_{i-1} 
\\ \nabla f_{2,x} \cdot \dot\gamma_{i-1} & \nabla
f_{2,y}\cdot \dot\gamma_{i-1} 
\end{pmatrix}
\dot\gamma_{i-1} \\ + 
(Df) \dot\gamma_{i-1} \times 
(Df)\ddot\gamma_{i-1} 
\end{equation}
where 
\begin{equation}\label{curv1a}
\nabla f_{1,x} = \begin{pmatrix} 
-2a + \varphi_{1, xx} \\ \varphi_{1, xy} \end{pmatrix}, 
\end{equation}
and 
\begin{equation}\label{curv1b}
\nabla f_{1,y} = \begin{pmatrix} 
\varphi_{1, xy} \\ \varphi_{1, yy} \end{pmatrix}; 
\nabla f_{2,x} = \begin{pmatrix} 
\varphi_{2, xx} \\ \varphi_{2, xy} \end{pmatrix}; 
\nabla f_{2,y} = \begin{pmatrix} 
\varphi_{2, xy} \\ \varphi_{2, yy} \end{pmatrix}; 
\end{equation}
We shall estimate the two terms of \eqref{curv1} separately. These
will yield the two terms in the statement of the Proposition. 
For the second term we have 
\[ 
|(Df) \dot\gamma_{i-1} \times 
(Df)\ddot\gamma_{i-1}| =
|Det (Df)| |\dot\gamma_{i-1} \times \ddot\gamma_{i-1}|=
|Det (Df)| \kappa_{i-1}|\gamma_{i-1}|^{3}.
\]
Indeed, for the first equality, 
\( |\dot\gamma_{i-1} \times \ddot\gamma_{i-1}| \) is the area of
the parallelogram defined by the two vectors \( \dot\gamma_{i-1}  \)
and \(  \ddot\gamma_{i-1} \), and 
\( |(Df) \dot\gamma_{i-1} \times 
(Df)\ddot\gamma_{i-1}| \) is the
are of the parallelogram defined by the two vectors 
\( \dot\gamma_{i-1} \) and \(\ddot\gamma_{i-1} \) which
of course just the image of the first parallelogram under \( Df \).
The second equality just follows immediately from the definition of \( 
\kappa_{i-1} \). So it just remains to show that 
the value of \( |Det (Df)|  \) is bounded above by some multiple of \( 
b \) and \( \eta \). Indeed, writing \( f=h+\varphi \) we have, by the
``row-linearity'' of the determinant,
\begin{align*}
Det(Df) &= Det \begin{pmatrix}
h_{1x}+\varphi_{1x} & h_{1y}+\varphi_{1y}  \\ 
h_{2x}+\varphi_{2x} &  h_{2x}+\varphi_{2y}  
\end{pmatrix} 
\\ &= 
Det \begin{pmatrix}
h_{1x}& h_{1y} \\ 
h_{2x}+\varphi_{2x} &  h_{2x}+\varphi_{2y}  
\end{pmatrix}  + 
Det \begin{pmatrix}
\varphi_{1x} & \varphi_{1y}  \\ 
h_{2x}+\varphi_{2x} &  h_{2x}+\varphi_{2y}  
\end{pmatrix} 
\\
& 
= Det \begin{pmatrix}
h_{1x}& h_{1y} \\ 
h_{2x}&  h_{2x} 
\end{pmatrix}  + 
Det \begin{pmatrix}
h_{1x}& h_{1y} \\ 
\varphi_{2x} &  \varphi_{2y}  
\end{pmatrix}  
+
Det \begin{pmatrix}
\varphi_{1x} & \varphi_{1y}  \\ 
h_{2x}&  h_{2x} 
\end{pmatrix} 
+
Det \begin{pmatrix}
\varphi_{1x} & \varphi_{1y}  \\ 
\varphi_{2x} &  \varphi_{2y}  
\end{pmatrix} 
\end{align*}
Using \( h_{1x}=-2a \), \( h_{1y}= 1 \), \( h_{2x}=b \), \( h_{2y}=0 \)
and \( \|\varphi\|_{C^{2}}\leq \eta \) this gives 
\[ 
Det (Df) \leq b+ (2a\eta + \eta) + (2a\eta + \eta) + \eta = 
b+ 4a \eta + 3 \eta \leq b + 12\eta
\]
where in the last step we have used the fact that \( a \) is close to \( 
2 \). Substituting this above gives the required bound for the second 
term of \eqref{curv1}. To bound the first term we write 
\[ 
(Df)\dot\gamma_{i-1} \times 
\begin{pmatrix} 
\nabla f_{1,x} \cdot \dot\gamma_{i-1} & \nabla f_{1,y} \cdot 
\dot\gamma_{i-1} 
\\ \nabla f_{2,x} \cdot \dot\gamma_{i-1} & \nabla
f_{2,y}\cdot \dot\gamma_{i-1} 
\end{pmatrix}
\dot\gamma_{i-1} 
= \begin{pmatrix}
a_{1} & b_{1} \\ c_{1} & d_{1}
\end{pmatrix}
\begin{pmatrix}
    v_{1} \\ v_{2} 
    \end{pmatrix} 
    \times 
    \begin{pmatrix}
    a_{2} & b_{2} \\ c_{2} & d_{2}
    \end{pmatrix}
    \begin{pmatrix}
v_{1} \\ v_{2} 
\end{pmatrix} 
\]
Then the norm of this expression is bounded above by 
\begin{align*}
 & |a_{1}c_{1}v_{1}^{2}+a_{1}d_{2}v_{1}v_{2} + b_{1}c_{2}v_{1}v_{2} + 
 b_{1}d_{2}v_{2}^{2} - a_{2}c_{1}v_{1}^{2} - a_{2}d_{1}v_{1}v_{2}
 - b_{2}c_{1}v_{1}v_{2} - d_{1}b_{2}v_{2}^{2}| 
 \\ \leq & 
 \max\{|a_{1}c_{2}-a_{2}c_{1}|, |b_{1}d_{2}-d_{1}b_{2}| +
 |a_{1}d_{2}-c_{1}b_{2}|+|b_{1}c_{2}-a_{2}d_{1}|\}
 (v_{1}^{2}+v_{2}^{2})
 \\ \leq &
 4 \max\{|a_{1}c_{2}-a_{2}c_{1}|, |b_{1}d_{2}-d_{1}b_{2}|, 
 |a_{1}d_{2}-c_{1}b_{2}|, |b_{1}c_{2}-a_{2}d_{1}|\}
 (v_{1}^{2}+v_{2}^{2})
 \\ \leq &
 8 \max\{|a_{1}c_{2}|, |a_{2}c_{1}|, |b_{1}d_{2}|, |d_{1}b_{2}|, 
 |a_{1}d_{2}|, |c_{1}b_{2}|, |b_{1}c_{2}|, |a_{2}d_{1}|\}
 |\dot\gamma_{i-1}|^{2}. 
 \end{align*}   
 All the terms contain a factor \( \dot\gamma_{i-1} \);
 each of the terms \(  b_{2}, c_{2}, d_{2}\), see \eqref{curv1b}, contains a bounded
 constant multiplied by the factor \( \eta \); the term \( a_{2} \),
 see \eqref{curv1a}, is of the order of \( 2a \) but here it is
 multiplied by either \( c_{1} \) or \( d_{1} \) each one of which
 contains a term which is bounded by \( \eta \). Therefore, there
 exists a constant \( K>0 \) such that 
 \[ 
 \left| 
 (Df)\dot\gamma_{i-1} \times 
 \begin{pmatrix} 
 \nabla f_{1,x} \cdot \dot\gamma_{i-1} & \nabla f_{1,y} \cdot 
 \dot\gamma_{i-1} 
 \\ \nabla f_{2,x} \cdot \dot\gamma_{i-1} & \nabla
 f_{2,y}\cdot \dot\gamma_{i-1} 
 \end{pmatrix}
 \dot\gamma_{i-1} 
 \right|
 \leq K\eta |\gamma_{i-1}|^{3}.
 \]
\end{proof}

\begin{proof}[Proof of Proposition \ref{smallcurv}]  
Applying Lemma \ref{curv} recursively we get 
\begin{align*}
\kappa_n(s) &\le  
K(b+\eta) \kappa_{n-1}(s) \frac{|\dot\gamma_{n-1}|^3 }{|\dot\gamma_n|^{3}}
   + K(b+\eta) \frac{|\dot\gamma_{n-1}|^3 }{|\dot\gamma_n|^{3}} 
   \\ &
   \le (K(b+\eta))^2 \kappa_{n-2}\frac{|\dot\gamma_{n-2}|^3 }{|\dot\gamma_n|^{3}} + 
   (K(b+\eta))^2\frac{|\dot\gamma_{n-2}|^3 }{|\dot\gamma_n|^{3}}
+K(b+\eta) \frac{|\dot\gamma_{n-1}|^3 }{|\dot\gamma_n|^{3}}
\\ &
 \le \ldots 
\end{align*}
Using the expansivity assumption and \( b, \eta \) small, this gives
\[
\kappa_n(s)  
\le \frac{1}{C^3} (K(b+\eta) e^{-\lambda})^n \kappa_0(s) + 
\frac{1}{C^3} \frac{K(b+\eta)
e^{-\lambda}}{1-K(b+\eta)e^{-\lambda}}\leq \kappa_{0}(s) \leq \alpha. 
\]
\end{proof}

\subsection{Critical points}
\label{hypcrit}

The next Proposition makes precise the notion of a 
   \emph{critical point of order }\( 
   k \).
   We recall that \( \gamma \) is a \( C^{2} \) \emph{admissible curve} if all
   its tangent vectors have slope \( <\alpha \) and it has curvature \(
   <\alpha \). 
   We say that \( \gamma \)
      is a \emph{long admissible curve} if it is an admissible curve which
      crosses the entire length of \( \Delta_{\varepsilon} \).
   
   \begin{proposition}\label{tang}
       Let \( \gamma\subset\Delta_{\varepsilon}\cap\mathcal D \) be a long admissible 
       curve. Then there exists a unique point \( c^{(k)}\in \gamma \) such
       that  \( \gamma_{0}= f(\gamma) \) has a 
       (quadratic) tangency at \( c^{(k)}_{0}=f(c^{(k)}) \in \mathcal
       V_{k}^{-}\cup\mathcal V^{+} \) 
       with the stable foliation 
       \(  \mathcal E^{(k)}\), for any \( k\geq k_{0} \).  
   Moreover there
   exists  a constant \( K \), independent of \( b, \eta \), such that 
   \( d(c^{(k)}_{0}, c_{0}^{(k+1)}) \leq Kb^{k}  \). In particular, 
   the sequence \( \{c^{(k)}_{0}\} \) is Cauchy. 
    \end{proposition}
    
    \begin{definition} 
	We call \( c^{(k)} \) and \( c_{0}^{(k)} \) respectively a \emph{critical
	point} and \emph{critical value} of order \( k \), associated to the long
	admissible curve \( \gamma \). 
     \end{definition} 

 We remark that critical values \( c_{0}^{(k)} \) of finite order are 
 not guaranteed to be outside \( \mathcal D \), however we shall show 
 below that their limit points as \( k\to\infty \), i.e. the ``real'' 
 critical points always fall strictly outside \( \mathcal D \) for \( 
 a> a^{*} \). 

Given a parametrized curve \( \gamma_{0}=\gamma_{0}(t) \) and its image \(
\gamma_{1}=\gamma_{1}(t) = f(\gamma_{0}(t)) \) we denote by \( \kappa_{0}(t) \)
the curvature of \( \gamma_{0} \) at the point \( \gamma_{0}(t) \) and
by \( \kappa_{1}(t)\) the curvature of \( \gamma_{1} \) at the point \( 
\gamma_{1}(t) \).

\begin{lemma}\label{poscurv}
   Let \( \gamma_{0}(t) \) be an admissible curve and let \(
   \gamma_{1}(t)=f(\gamma(t)) 
   = (\xi_{1} (t), \eta_{1}(t)) \). Suppose that for some \( t \) we
   have \( \dot \eta_{1}(t) \neq 0 \) and 
   \( |\dot \xi _{1}(t)/ \dot \eta_{1}(t)| < 1 \).  Then \( |\kappa_{1}(t)| > 
   a/b \gg 1 \). 
 \end{lemma}   

 Lemma \ref{poscurv} essentially says that if the tangent direction 
of the image of an admissible curve at a certain point is roughly
vertical (or at least contained in the ``vertical'' cone between the
positive and the negative diagonals) then the curvature at this point 
is strictly bounded away from 0.  This does not apply to admissible
curves outside \( \Delta_{\varepsilon} \) since we have shown above (
Corollary \ref{admisstoadmiss}) that images of such curves are still
admissible and therefore their tangent directions are roughly
horizontal. We will instead apply it below to the images of admissible curves
inside \( \Delta_{\varepsilon} \) as a way of pinpointing the location
of \emph{folds}.

\begin{proof}
First recall that the curvature \( \kappa_{1}(t) \) is independent of the
choice of parametrization and also the condition 
   \( |\dot \xi _{1}(t)/ \dot \eta_{1}(t)| < 1 \)  is independent of 
   the parametrization since \( |\dot \xi _{1}(t)/ \dot \eta_{1}(t)| \)
   is just the slope of the tangent vector.
Therefore we choose the parametrization 
\[ 
\gamma_{0}(t) = (t, y(t)). 
\] 
For simplicity we also omit the
subscript \( 1 \) from the coordinate functions of \( \gamma_{1} \)
and just write \( \gamma_{1}(t) = (\xi(t), \eta(t)) \).
From the definition of \( f \) we have 
\begin{align*}
    (\xi(t), \eta(t)) &= (1+at^{2}+y(t) + \varphi_{1}(\gamma_{0}(t)),
    bt+\varphi_{2}(\gamma_{0}(t)))
    \\
    (\dot\xi(t), \dot\eta(t)) &= (-2at + \dot y(t) +
    \nabla\varphi_{1}(\gamma_{0}(t)) \cdot \dot\gamma_{0}(t), 
    b+\nabla\varphi_{2}(\gamma_{0}(t)) \cdot \dot\gamma_{0}(t))
    \\
    (\ddot \xi(t), \ddot\eta(t)) &= 
    (-2a + \dot y(t) + D^{2}\varphi_{1}(\gamma_{0}(t)) [\dot\gamma_{0}(t)]^{2}, 
    D^{2}\varphi_{2}(\gamma_{0}(t))[\dot\gamma_{0}(t)]^{2})
    \end{align*}
  Choosing \( \eta \) sufficiently small, for example so that \(
  4\|\nabla\varphi_{2}(\gamma_{0}(t))\|(1+\alpha) < b \)  this implies
  \begin{equation}\label{etadot}
      3b/4 \leq |\dot\eta(t)| \leq 5b/4.
      \end{equation}
We can now compute the curvature \( \kappa_{1}(t) \). First notice
that the condition   \( |\dot \xi _{1}(t)/ \dot \eta_{1}(t)| < 1 \)
implies in particular \( \|(\dot\xi(t), \dot\eta(t))\|\leq \sqrt 2
|\dot\eta(t)| \). Then we have 
\[ 
\kappa_{1}(t) = 
\frac{|\ddot \xi (t) \dot\eta (t) - \dot\xi(t) \ddot\eta
(t)|}{\|(\dot\xi(t), \dot(\eta(t))\|^{3}}\geq 
\frac{|\ddot \xi (t) \dot\eta (t) - \dot\xi(t) \ddot\eta
(t)|}{4 |(\dot\eta(t))|^{3}}
\]
Dividing numerator and denominator by \( |\dot\eta(t)| \), using
the condition \( |\dot \xi _{1}(t)/ \dot \eta_{1}(t)| < 1 \) and
\eqref{etadot} we get 
\[ 
\kappa_{1}(t) \geq 
\frac{|\ddot \xi (t)  - \frac{\dot\xi(t)}{\dot\eta(t)} \ddot\eta
(t)|}{4 |(\dot\eta(t))|^{2}}
\geq 
\frac{|\ddot \xi (t)| - \left|\frac{\dot\xi(t)}{\dot\eta(t)}\right| \  
|\ddot\eta (t)|}{4 |(\dot\eta(t))|^{2}}
\geq 
\frac{|\ddot \xi (t)| - |\ddot\eta (t)|}{4 |(\dot\eta(t))|^{2}}
\geq 
\frac{|\ddot \xi (t)| - |\ddot\eta (t)|}{7b^{2}}
\]
  Finally, from the formulas for \( \ddot\xi(t) \) and \( \ddot\eta(t) \)
  and the fact that \( |\dot y(t)| \leq \alpha  \) by the
  admissibility of \( \gamma_{0} \), 
  we get 
  \[ 
  |\ddot \xi (t)| - |\ddot\eta (t)| \geq 2a - \alpha - 2
  \|\varphi\|_{C^{2}} \geq a
  \]
  as long as \( \eta \) is sufficiently small. 
      
\end{proof}

   \begin{proof}[Proof of Proposition \ref{tang}]

The existence of a tangency between \( f(\gamma) \) and the stable
foliation \( \mathcal E^{(k)} \) 
 follows by the simple geometric observation that the image 
of a  long admissible curve necessarily ``changes direction'' between 
one endpoint and the other. Thus,  
by a simple Intermediate Value argument it follows that
there is some point of tangency. 

Now, Proposition
\ref{c2close}  says that the leaves of the stable foliations \(
\mathcal E^{(k)} \) are close to the piece of stable manifold 
\( f^{-1}(W^{s}_{\delta}(q) \)  and thus have slope close to \( 2 \),
and that
their curvature is small. In
particular the point of tangency must occur at some point at which the
tangent direction to \( f(\gamma) \) is close to 2 and therefore
Proposition  \ref{tang}  shows that at this point of tangency 
\( f(\gamma) \) has strictly positive curvature.  This implies
that this tangency is quadratic as well as unique.

  \end{proof}

\begin{figure}
    \includegraphics[width=8cm]{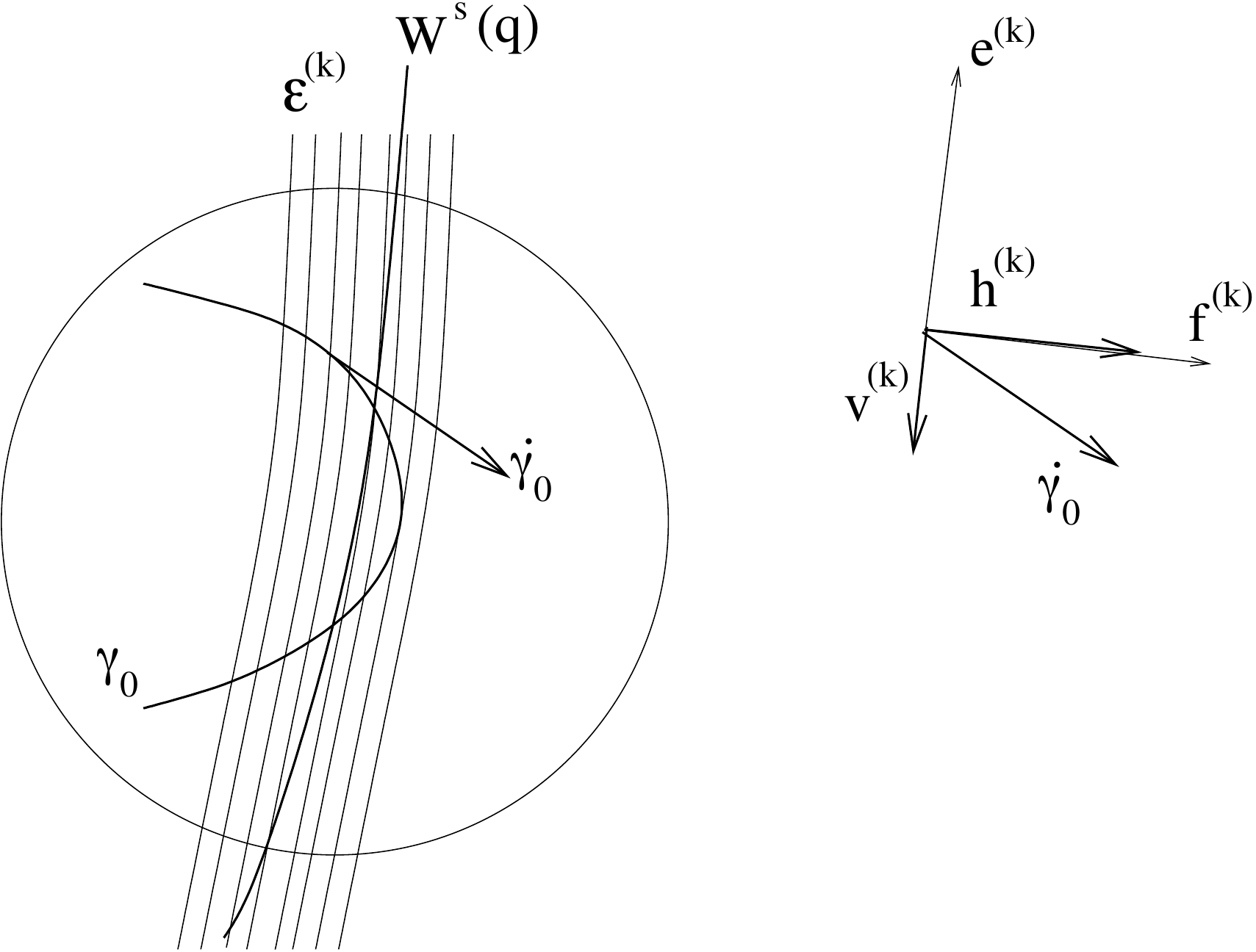}
    \caption{Hyperbolic coordinates}
     \end{figure}

\section{Hyperbolicity estimates}
\label{sectionhypest}

This is the final and main section of the paper. We apply the notion
of hyperbolic coordinates and dynamical defined critical points to
prove Theorem 1. 
In section \ref{shadowing} we
combine the hyperbolic coordinates and the curvature estimates to show
that all components of the unstable manifold \( W^{u}(p) \) 
in \( \Delta_{\varepsilon} \) 
are almost horizontal curves with
small curvature. In particular they all have well-defined critical
points.  In section \ref{HypRet} we take advantage of the structure of
critical points on such components to show that points in the critical
region \( \Delta_{\varepsilon}\setminus\Delta \) recover hyperbolicity 
after some bounded number of iterations depending only on the
parameter \( a \). In Section \ref{UnifHyp} we then extend these
estimates to uniform expansion estimates on all of \( W^{u}(p) \) with
a hyperbolicity constant \( C_{a} \) depending only on the parameter. 
In Section \ref{UnifHypOm} we then show how to extend this
hyperbolicity to the closure of \( W^{u}(p) \) and thus to the whole
nonwandering set \( \Omega \). Finally, in Section \ref{Lyap}, we
consider the bifurcation parameter value \( a=a^{*} \) and show that
all Lyapunov exponents are uniformly bounded away from 0.

\subsection{Shadowing}
\label{shadowing}

Let 
\begin{equation}\label{lambda}
\lambda = \min\left\{\frac{1}{2}\ln \frac{3}{\sqrt 5}, \hat\lambda\right\}.
\end{equation}

\begin{proposition}   \label{recovering}
For all \( a\geq a^{*} \) all components of \(
W^{u}(p) \cap \Delta_{\varepsilon}  \) are long admissible curves. 
Moreover, for all \( z\in
W^{u}(p)\cap(\Delta_{\varepsilon}\setminus\Delta) \) and any
vector \( v \) tangent to \( W^{u}(p) \) at \( z \) and \( k\geq 1 \) 
such that \( f(z)\in \mathcal V_{k}\setminus \mathcal V_{k+1} \) we have 
\[ 
\|Df^{k}_{z}(v)\|\geq e^{\lambda k}\|v\|. 
\] 
\end{proposition}
We emphasize that Proposition \ref{recovering} holds 
also for parameter values for which the first tangency occurs. 

\begin{proof}
We  first prove the expansivity statement and then the admissibility
of leaves of \( W^{u}(p) \) in \( \Delta_{\varepsilon} \). 

\subsubsection{Expansion}
If \( \gamma(s) = (x(s), y(s)) \subset\Delta_{\varepsilon}\cap \mathcal D \) 
is a long
admissible curve we consider the tangent vectors \( \dot\gamma(s) \)
and their images \( \dot\gamma_{0}(s) = Df(\dot\gamma(s)) \). By 
Proposition \ref{tang}, \( \dot\gamma_{0} \) is tangent 
to the stable direction \( e^{(k)} \) at the point \( c_{0}^{(k)} \). 
For this and other nearby points on \( \gamma \) we can write the
tangent vector 
\[ 
\dot\gamma_{0} = h^{(k)}_{0}f^{(k)} + v^{(k)}_{0} e^{(k)}
\]
where \( (f^{(k)}, e^{(k)}) \) is the orthogonal basis given by the
most expanded and most contracted direction for \( Df^{k} \) and \(
h_{0}^{(k)} \) and \( v_{0}^{(k)} \) are the components of \(
\dot\gamma_{0} \) in this basis. Notice that the basis itself depends 
on the point. Proposition \ref{c2close} implies that the basis
varies \emph{very slowly} with the base point, and 
Proposition
\ref{tang} implies that the tangent vector \( \dot\gamma_{0} \) is
varying at \emph{positive speed} with respect to this basis. We omit
the calculations which are relatively standard,   see for example
\cite{LuzVia03}.
Specifically this implies that the component \( h_{0}^{(k)} \)
of the tangent vector \( \dot\gamma_{0} \) 
at some point \( z_{0}=f(z) \in \gamma_{0} \)
is proportional to the distance between \( z \) and the critical point
of order \( k \), \( c^{(k)} \).  In our setting, the constants
actually give 
\begin{equation}\label{horizontal}
    |h_{0}^{(k)}(z_{0})| \geq d(z, c^{(k)}).
 \end{equation}   
We can now prove the following
\begin{lemma}\label{binding}
Suppose \( \gamma \subset \Delta_{\varepsilon} \) is an admissible
curve, \( z\in \gamma \), \( z_{0}=f(z)\in\mathcal
V_{k}\setminus\mathcal V_{k+1} \) and \( c^{(k)} \) is the critical
point of order \( k \) in \( \gamma \). 
Then for a vector \( w \) tangent to \( \gamma \) at \( z \) and all \( 
j=0,\ldots, k \) we have 
\[ 
\|Df^{j+1}_{z}(w)\| \geq 3^{j} d(z, c^{(k)})\|w\|.
\]
In particular
\[ 
\|Df^{k+1}_{z}(w)\|\geq e^{\lambda (k+1)}\|w\|.
\]
\end{lemma}
\begin{proof}
  The first equality follows immediately from   \eqref{horizontal}
  and \eqref{contest}.  To prove the second we need to find a bound for
  \( d(z, c^{(k)}) \) in terms of \( k \). Using the quadratic nature 
  of \( \gamma_{0} \) and the proximity to the one-dimensional map 
  \( 1-ax^{2} \) with \( a\approx 2 \), we  obtain 
  \begin{equation}\label{critdist}
  d(z, c^{(k)}) \geq \frac{1}{3} \sqrt{d(z_{0}, c_{0}^{(k)})}.
  \end{equation}
  To estimate \( d(z_{0}, c_{0}^{(k)}) \) we use the observation that 
  the ``real'' critical value \( c_{0} \) on \( \gamma_{0} \), i.e.
  the point of tangency between \( \gamma_{0} \) and the limiting
  stable foliation \( \mathcal E^{(\infty)} \) lies necessarily either on \(
  W^{s}(q) \) (this is only a possibility if \( a=a^{*} \)) or to the 
  right of \( W^{s}(q) \) in \( \mathcal Q \). We write this as \(
  \delta_{0}= d(c_{0}, W^{s}(q)) \geq 0 \). Combining this with 
  Lemma \ref{distance1} and the rate of convergence 
  of critical points of finite order \( d(c_{0}^{k}, c_{0}) \leq K b^{k} \)
  as mentioned in Proposition  \ref{tang} and taking \( b \)
 sufficiently small, we get 
\begin{align*}
d(z_{0}, c_{0}^{(k)}) &\geq d(z_{0}, W^{s}(q)) + d(W^{s}(q), c_{0}) - 
d(c_{0}^{(k)}, c_{0}) \\ &\geq \frac{\delta}{2}5^{-k}+\delta_{0}- K b^{k} 
\geq \frac{\delta}{3}5^{-k}.
\end{align*}
Substituting this into \eqref{critdist}  and using the fact that we
can assume \(
k\geq k_{0} \) as well as the definition of \( k_{0} \) in \eqref{k0}
and of \( \lambda \) in \eqref{lambda},
we have 
\[ 
3^{k}d(z,c^{(k)}) \geq \frac{\sqrt \delta}{2\sqrt
3}\left(\frac{3}{\sqrt{5}}\right)^{k}\geq e^{\lambda (k+1)}
\]
\end{proof}  

\subsubsection{Admissibility}
Returning to the proof of the Proposition, 
to obtain the statement about admissibility, notice first of all that 
combining Lemma \ref{binding} with Lemma \ref{smallcurv} we immediately
obtain the statement that 
 if \( \gamma\subset W^{u}(p)\cap\Delta_{\varepsilon} \) is admissible
 and \( k \) is the first time that \( f^{k}(\gamma) \subset
 \Delta_{\varepsilon} \), then \( f^{k}(\gamma) \) is admissible. 
Now, by choosing \( |b| \) and \( \eta \) small we can guarantee that  \( 
W^{u}_{loc}(p)\cap\Delta_{\varepsilon} \) is a long admissible 
curve. Moreover, every piece of \( W^{u}(p) \cap\Delta_{\varepsilon} \)  
is the image of some curve in \(    W^{u}_{loc}(p)\cap\Delta_{\varepsilon}  \)
and is therefore admissible. 
\end{proof}

\subsection{Hyperbolicity after returns to  \protect\(
\Delta_{\varepsilon} \protect \)}
 \label{HypRet}
 
 Proposition \ref{recovering} gives a pointwise recovery time for 
 the hyperbolicity of points in the critical region, based on their
 position. The following Proposition gives a key uniformity estimate 
 in the phase space for each parameter \( a> a^{*} \). 
 
\begin{proposition} \label{N}
For all \( a> a^{*} \) there exists a constant \( N_{a} 
\) such that for \( z\in W^{u}(p)\cap\Delta_{\varepsilon}\cap 
\Omega(f)\), and \( v \) a tangent vector to \( W^{u}(p) \) at \( z
\), there exists \( n(z) \leq N_{a} \) such that \( Df^{n(z)}_{z}(v) \) is
almost horizontal and 
   \[ 
   \|Df^{n(z)}_{z}(v)\|\geq e^{\lambda n(z)}\|v\|. 
   \] 
\end{proposition}

We remark that the constant \( N \) is \emph{not} 
uniformly bounded in \( a\) and in particular 
\emph{does not} apply to
\( a=a^{*} \).  However  
it gives us a \emph{uniformity statement} in \( z \) which will
implies, as we shall see below, 
uniform hyperbolicity for each given parameter value \( a> a^{*} \).
For the proof we need to
extend the definition of admissibility naturally to curves which
are only differentiable of class \( C^{1+1} \) (Lipschitz continuous
derivative). 

\begin{definition}
We say that \( \gamma(s)\subset \Delta_{\varepsilon} \) is  a \( C^{1+1} \) 
admissible curve if \( |\dot y| /|\dot x| <
    \alpha  \), and \( \dot \gamma (s) \) is Lipschitz with Lipschitz 
    constant \( \leq  \alpha \). 
\end{definition}

We also give the formal definition of a ``real'' critical point, which
applies both to \( C^{2} \) and to  \( C^{1+1} \) admissible
  curves.
  
  \begin{definition}
  We say that \( c\in\gamma \) is a critical point if \( e^{(\infty)} \)
    is defined at \( f(c)\in \gamma \) and coincides with \(
    Df_{c}(\dot\gamma(c)) \). 
\end{definition}

    \begin{lemma}\label{unique critical}
For every \( a>a^{*} \), 
every \( z\in \overline{W^{u}(p)}\cap\Delta_{\varepsilon}\cap\Omega\) 
lies on a \( C^{1+1} \) admissible curve \( \gamma \) which is the 
limit of \( C^{2} \) admissible curves in \( W^{u}(p) \) and \( \gamma \)
contains a unique critical point \( c(\gamma) \) with 
\( d(z, c)>0 \). 
\end{lemma}

\begin{proof} We split the proof into two parts. 
    
    \subsubsection*{Every point lies on an admissible curve}
    We show first of all that every point 
    \( z\in \overline{W^{u}(p)}\cap\Delta_{\varepsilon}\cap\Omega \) 
lies on a \( C^{1+1} \) admissible curve which is the limit of \( C^{2} \)
admissible curves in \( W^{u}(p) \). 
    Let \( z\in \overline{W^{u}(p)}\cap\Delta_{\varepsilon}\cap\Omega \) 
    and let \( z_{n}\to z \) be a sequence in 
    with \( z_{n}\in W^{u}(p)\cap\Delta_{\varepsilon}\cap\Omega \). 
    By Proposition \ref{recovering} each \( z_{n} \) belongs to 
    a long admissible curve \( \gamma_{n}\subset W^{u}(p) \). We can
    write these as functions \( 
    \gamma_{n}:I\to \mathbb R \) with \( I=[-\varepsilon, \varepsilon] \)
    and suppose that converge pointwise to \( 
    \gamma: I \to \mathbb R\). Since \( I \) is compact and \( 
    \gamma_{n}, \dot\gamma_{n} \) are bounded and equicontinuous 
    sequences we have that \( \gamma \) is \( C^{1} \) and \( 
    \gamma_{n}\to \gamma \) in the \( C^{1} \) topology. 
    To see that \( \dot\gamma \) is Lipschitz, let \( x,y \in I \) and 
    observe that each \( \dot\gamma_{n} \) is  a Lipschitz function 
    with uniformly bounded Lipschitz constant \( \alpha
    \). Then we have \( |\dot\gamma_{n}(x)-\dot\gamma_{n}(y)|\leq 
    \alpha |x-y| \) and hence \( |\dot\gamma_{n}(x)-\dot\gamma_{n}(y)|\leq 
    \alpha |x-y| \). 

    \subsubsection*{Every admissible curve contains a critical point}
    We now show that any such curve \( \gamma \) contains a unique
    critical point. We show first that it must contain at most one,
    and then argue that it must contain at least one. 
Let \( \theta (\gamma_{n}(t)) \) be the angle between the vectors \( 
     Df_{(t, \gamma_{n}(t))}(1, \gamma'_{n}(t)) \) and \( e^{\infty}(f(t, 
     \gamma_{n}(t)) \). Since the image of each admissible curve is 
     quadratic with respect to \( \mathcal E^{(\infty)} \) we have that 
     \( \theta (\gamma_{n}(t)) \) has a strictly non-zero derivative 
     having at most one root corresponding to a point of tangency 
      between \( f(\gamma_{n}) \) and \( \mathcal E^{(\infty)} \).  
      Since \( \gamma_{n}\to\gamma \) in the \( 
     C^{1} \) topology, we have that \( \theta(\gamma(t)) \) also has 
     strictly non-zero derivative having at most one root also 
     corresponding to a point of tangency 
      between \( f(\gamma) \) and \( \mathcal E^{(\infty)} \).
     To see that such a point
     exists, observe that if \( a> a^{*} \) then the 
     graph of \( 
     \gamma \) crosses the boundary of \( \Delta \) twice 
     and \( f(\gamma\cap\Delta) \) is outside \( \mathcal D \) 
     where the foliation \( \mathcal E^{(\infty)} \) is well defined 
     and the extreme points of \( f(\gamma\cap\Delta) \) both lie on 
     a piece of \( W^{s}(q) \) which is a leaf of the foliation \( 
     \mathcal E^{(\infty)} \). This implies that there exists a point 
     outside the interior of \( \mathcal D \) 
     where \( f(\gamma) \) is tangent to \( 
     \mathcal E^{(\infty)} \).  
\end{proof}

Lemma \ref{unique critical} allows us to define a canonical set \(
\mathcal C_{a} \) of \emph{critical points} as the union of all critical
points \( c(\gamma) \) for every \( C^{1+1} \) which are \( C^{1} \)
limits of 
long admissible curves of \( W^{u}\cap \Delta_{\varepsilon} \). 
In the next Lemma we show that this set is bounded away from the set
of nonwandering points.

\begin{lemma}\label{crit}
For all \( a> a^{*} \) we have \( d(\mathcal C_{a}, \Omega) > 0 \).
\end{lemma}

We emphasize that \( d(\mathcal C_{a}, \Omega) \) is not uniformly bounded 
in the parameter. The constant \( N_{a} \) in the Proposition will be
defined below in terms of \( d(\mathcal C_{a}, \Omega) \). 

 \begin{proof}
     Notice first of all that \( \mathcal C_{a}\subset\Delta_{\varepsilon} 
      \) and thus in particular is bounded. Let \( c_{k}=c(\gamma_{k}) \) 
      be a sequence converging to some point \( c \). We need to show 
      that \( c\in\mathcal C_{a} \). Since each \( \gamma_{k} \) is the 
      limit of long admissible curves, we can consider sequences \( 
      \gamma_{k}^{(n)}\to \gamma_{k} \) for each \( k \). 
      Using Lemma \ref{unique critical} and the fact that \( 
      \{\gamma_{k}^{(k)}\} \) converges pointwise to \( \gamma \), we 
      conclude that this convergence is in fact \( C^{1} \). Since 
      \( \theta(\gamma_{k}^{(k)}(c_{k}))\to 0 \) we have that \( 
      \theta (\gamma(c))=0 \) and this implies that \( c \) is a 
      critical point as required. 
      
      We have therefore shown that the critical set \( \mathcal C_{a}     \) 
      is compact. Since \( \Omega \) is also compact, it is sufficient
      to show that \( \mathcal C_{a}\cap\Omega = \emptyset \) to imply
      that they are at some positive distance apart. Disjointness
      follows from the observation that 
the image of a critical point is 
always outside \( \mathcal D \), while \( \Omega \) is an 
invariant set contained in \(  \mathcal D \).


  \end{proof}

\begin{proof}[Proof of Proposition \ref{N}]
 By Lemma \ref{crit} and the uniform approximation of the critical set
 \( \mathcal C \) by the finite order critical sets \( \mathcal C^{(n)}
 \), 
 there exists \( N_{a} \) sufficiently large so that the following
 two conditions hold (using also \( \lambda < \log 3) \): 
 \begin{equation}\label{Ndef}
     d(\mathcal C^{(N_{a})}_{a}, \mathcal C_{a}) < d (\mathcal C_{a}, 
     \Omega) /2 \  \ 
     \text{ and } \  \ 
     3^{N_{a}} d(\mathcal C^{(N_{a})}_{a}, \Omega)
     \geq e^{\lambda N_{a}}.
    \end{equation}
 Now consider 
\( z\in \Delta_{\varepsilon}\cap W^{u}(p)\cap
    \Omega \) and let \( n\geq 1 \) be 
    such that \( f(z)\in \mathcal
    V_{n}\setminus \mathcal V_{n+1} \).  
Recall \( f(\Delta_{\varepsilon}) \subset \mathcal
   V_{k_{0}}\) and therefore such an \( n \) is always well defined
   except for those points which map exactly to the the curve \(
   f^{-1}(W^{s}_{\delta}(q) \) which forms the boundary between \(
   \mathcal V^{+} \) and \( \mathcal V^{-1} \). For these points we
   let \( n=+\infty \). 
   Then we let 
    \[ 
    n(z) = \min\{n, N_{a}\}.
    \]
    If \( n \leq N_{a} \) the statement follows from Proposition
    \ref{recovering}. 
   Otherwise our choice of \( N_{a} \) in \eqref{Ndef}
    gives 
    \[ 
    \|Df^{N_{a}}(v)\| \geq  
    3^{N_{a}}d (z, \mathcal C_{a}^{(N_{a})})\|v\|
    \geq 
    3^{N_{a}}d (\Omega, \mathcal C_{a}^{(N_{a})})\|v\|
    \geq 
    3^{N_{a}}d (\Omega, \mathcal C_{a})\|v\| /2
    \geq e^{\lambda N_{a}}\|v\|
    \]
The first inequality follows from Lemma \ref{binding}, the second one 
follows from \( z\in \Omega \), the third one follows from the first
part of \eqref{Ndef}, and the last one follows from the second part of
\eqref{Ndef}.
     
Finally, considering the components of 
\( v \) in hyperbolic coordinates we have 
\( \|v_{N_{a}}^{(N_{a})}\|\leq (b/3)^{N_{a}} \) 
     and \( \|h_{N_{a}}^{(N_{a})}\|\geq e^{\lambda N_{a}} \) and therefore \(
     Df^{N_{a}}(v) \) is almost horizontal. 

    \end{proof}
    
\subsection{Uniform hyperbolicity on \protect\( W^{u}(p) \protect\)}
    \label{UnifHyp}
    
 The following Proposition is is essentially a Corollary of
 Proposition \ref{N}. However we state it separately as it gives an
 explicit construction of the constant \( C_{a} \) of hyperbolicity
 for each \( a> a^{*} \).    Before stating the result we define this
 constant.

 Let 
 \(
 C^{-}_{N_{a}}=\min\{\|(Df^{j}_{z})^{-1}\|^{-1}: x\in\mathcal D, 1\leq j \leq 
 N_{a}\}
 \)
 and 
 \(
 C^{+}_{N_{a}}=\max\{\|Df^{j}_{z} \|: x\in\mathcal D, 1\leq j \leq 
 N_{a}\}
 \)
 denote 
 the maximum possible contraction and the maximum possible expansion 
 exhibited by any vector \( 
 v\in T_{x}\mathbb R^{2} \) for any point \( x\in\mathcal D \) in at 
 most \( N_{a} \) iterations. Letting \( C_{\varepsilon} \) denote the 
 constant of hyperbolicity as in \eqref{UE1} on page \pageref{UE1}, 
 we then let 
 \[ 
 C_{a}=\min\left\{\frac{C_{\varepsilon}}{C^{+}_{N}},
 \frac{C_{N}^{-}e^{-\lambda N}}{C^{+}_{N}} \right\} 
 \]

\begin{proposition}\label{HypProp1}
    For all \( a> a^{*} \), all \( z\in W^{u}(p)\cap\Omega(f) \) 
    and all vectors \( w \) tangent to \( W^{u}(p) \) at \( z \)
    we have 
    \[ 
    \|Df^{n}_{z}(w) \| \geq  C_{a} e^{\lambda n}\|w\|
    \]
    for all \( n\geq 1 \). 
    \end{proposition}
    
\begin{proof} 
Let    \( z\in W^{u}(p)\cap\Omega(f) \)  and let 
\( w \) tangent to \( W^{u}(p) \) at \( z \). Since we do not assume
anything about the location of \( z \) the vector \( w \) may or may
not be almost horizontal. We distinguish these two possibilities.

 \subsubsection*{Case 1: \protect\( w \protect\) almost horizontal}
Let
 \( 0 \leq k_{1}< \ldots < k_{s} < n \) be the sequence of returns of 
 the iterates of \( z \) to \( \Delta_{\varepsilon} \) (
with \( k_{1}=0 \) if \( z\in\Delta_{\varepsilon } \) and \( k_{1}>0 \)
 otherwise). Then for each \( k_{i} \) we have an integer \(
 n_{i}=n(z_{k_{i}}) \leq N_{a}\) given by Proposition \ref{recovering}. 
 Then we can write 
 \[ 
 k_{i+1}=k_{i}+n_{i}+q_{i}
 \]
 where \( q_{i} \) is the number of iterates during which the point
 remains outside \( \Delta_{\varepsilon} \). From Proposition
 \ref{recovering} properties \eqref{slope} and \eqref{UE2}, the images
 of the vector at these iterates remains horizontal and we have 
 \[ 
 \|Df^{k_{i}}_{z}(w)\|\geq e^{\lambda k_{i}}\|w\|
 \]
 for all \( i=1,\ldots, s \), in particular for \( i=s \). 
 If \( k_{s}+n_{s}\leq n \), applying \eqref{UE2}  to the remaining
 iterates gives \( \|Df_{z}^{n}(w)\|\geq C_{\varepsilon}e^{\lambda
 n}\|w\| \geq C_{a}e^{\lambda
 n}\|w\| \) as required. 
 
 If \( k_{s}+n_{s} > n \) we have expansion for the first \( k_{s} \) 
 iterates which gives  
 \( \|Df^{k_{s}}(w)\| \geq e^{\lambda  k_{s}}\|w\| \). There follow 
 \( n-k_{s} \leq n_{s}\leq N_{a} \) iterates (since \( n_{s}\leq
N_{a} \)) during which we can bound the contraction
 coarsely by the \( N_{a} \)'th power of the maximum contraction 
 in the region \( \mathcal D \) which gives 
 \[ 
 \|Df^{n}(w)\|\geq C_{N}^{-} e^{\lambda k_{s}}\|w\| = C_{N}^{-}
 e^{-\lambda N} e^{\lambda n}\|w\|.
 \]
 
 \subsubsection*{Case 2: \protect\( w \protect\) is not almost horizontal}
 We now suppose that \( w \) is not almost horizontal.  
 \begin{claim}
     There exists 
     \[
     N_{a}\geq m > 0
     \]
     such that \( f^{-m}(z) \in \Delta_{\varepsilon} \)
and \( w_{-m}= Df^{-m}(w) \) is almost horizontal. 
\end{claim}
\begin{proof}
We show first of all  that some preimage of \( z \) lies in \(
 \Delta_{\varepsilon} \). Indeed, \( z\in W^{u}(p) \) implies that \(
 z_{-n}\to p \) as \( n\to \infty \) and therefore that \( w_{-n} \)
 is almost horizontal for sufficiently large \( n \) since the local
 unstable manifold of \( p \) is admissible.  By the invariance of the
 unstable conefield outside \( \Delta_{\varepsilon} \) images of \(
 w_{-n} \) remain almost horizontal unless some return to \(
 \Delta_{\varepsilon} \) takes place. 
 
 Now let \( m>0 \) be the
 smallest integer such that \( f^{-m}\in \Delta_{\varepsilon} \). Then
 \( w_{-m} \) is almost  horizontal since every component of \( W^{u} \)
 in \( \Delta_{\varepsilon} \) is almost horizontal. 
 From Proposition \ref{N} it follows that \(
  Df_{z_{-m}}^{n(z_{-m})}(w_{-m}) \) is almost horizontal and
  therefore it follows necessarily that \( m\leq n(z_{m}) \leq N_{a} \).
  Otherwise \( w\) will be almost horizontal.
  \end{proof}
 Returning to the proof of the Proposition, 
 we can now  argue as in the previous case to obtain
 exponential growth starting from time \( -m \): 
 \begin{equation}\label{C}
 \|Df^{n}(w)\|=\|Df^{n+m}(w_{-m})\| \geq C' e^{\lambda (n+m)}
 \|w_{-m}\|
 \end{equation}
 where \( C'= \min\{C_{\varepsilon}, C_{N}^{-} e^{-\lambda N}\}. \) 
Moreover
 \[
 \|w\|=\|Df^{m}(w_{-m})\|\leq \|Df^{m}\| \ \|w_{-m}\| \leq C_{N}^{+}
 \|w_{-m}\|.
 \]
Substituting this back into \eqref{C} completes the proof. 
\end{proof}     
    
\subsection{Uniform hyperbolicity on \protect\( \Omega \protect \)}
  \label{UnifHypOm}
  
 We have obtained  uniform expansion estimates for vectors
 tangent to \( W^{u}(p) \). In this section we show that these
 estimates can be extended to \( \Omega \). 
 This part of the argument uses very little of the 
 specific H\'enon-like form of the map and therefore we state it in a 
 more abstract and general context.

\begin{proposition}\label{HypProp2}
Let \( f: \mathbb R^{2}\to \mathbb R^{2} \) be a \( C^{1} \)
diffeomorphism and \( \Omega \) a compact invariant set with \( |\det 
Df| < 1 \) on \( \Omega \). Suppose that 
there exists  some invariant submanifold \( W \) dense in \( \Omega \)
and such that 
there exist constants \( C, \lambda > 0 \) such that \(
    \|Df_{z}(v)\|\geq Ce^{\lambda n} \) for all \( z\in W\cap \Omega \) and 
    \( v \) tangent to \( W \).
 Then \( \Omega \) is uniformly hyperbolic with hyperbolic constants \( 
 C \) and \( \lambda \).  
\end{proposition}    

Proposition \ref{HypProp2} completes the proof of part \( (a) \) of
the Theorem and shows that the rates of expansion and contraction 
admit uniform bounds independent of the parameter. 

\begin{proof}
We shall show that \( \Omega \) is uniformly hyperbolic by
constructing an invariant hyperbolic splitting \( E^{s}_{z}\oplus
E^{u}_{z} \) at every point of \( \Omega \) and then showing that this
splitting is continuous. We carry out this construction in several
steps. 
The starting point is the observation that \( E^{u}_{z} \) is
already given by the tangent direction to \( W \) for all points \(
z\in \Omega\cap W \). 

\begin{lemma}\label{limit}
    For any \( z\in\Omega \) and any sequence \( z_{j}\in W \) with \( 
    z_{j}\to z \), the sequence \( E^{u}(z_{j}) \) converges to a
    (unique) limit direction \( E^{u}(z) \). 
    Each vector \( v\in E^{u}(z) \) satisfies 
    \[ 
    \|Df_{z}^{n}(v) \| \geq C e^{\lambda n}\|v\| \quad \text{ and
    }\quad 
    \|Df_{z}^{-n}(v) \|\leq C^{-1}e^{-\lambda n}\|v\|
    \]
    for all \( n\geq 1 \). 
  \end{lemma}  
  \begin{proof}
Suppose \( z\in\Omega \) and let \( z_{j}\in W \) be a sequence
with \( z_{j} \to z \). Consider the sequence of corresponding
directions \( E^{u}(z_{j}) \). By compactness (of the space \( \mathbb 
S^{1} \) of possible directions) there must exist some subsequence \( 
z_{j_{i}} \) such that the corresponding directions \( E^{u}_{j_{i}} \)
converge to some direction which we call \( E^{u}(z) \). Notice that a
priori this direction is not unique since it depends on the choice of 
subsequence.  We shall show first that the forward expansion and
backward contraction estimates hold and then show that this actually
implies uniqueness. 

Let \( v\in E^{u}_{z} \) and \( v_{j_{i}}\in E^{u}_{z_{j_{i}}} \) be a
sequence with \( v_{j_{i}} \to v \). Then, for each \( 
    n\in\mathbb N \) we have, by the continuity of \( Df^{n} \), 
    \[ 
    \|Df^{n}_{z_{j}}(v_{j})\|\to \|Df^{n}_{z}(v)\| 
    \]
By assumption we know that \( \|Df_{z_{j_{i}}}(v_{j})\|\geq
Ce^{\lambda n}\|v_{j}\| \) and therefore 
    \[ 
    \|Df^{n}_{z}(v)\| \geq Ce^{\lambda n}-\varepsilon
    \]
    for any \( \varepsilon > 0 \). Therefore \( \|Df^{n}_{z}(v)\| 
    \geq Ce^{\lambda n} \) and, since this holds for every \( n \),  we 
    have the required statement as far as the expansion in forward
    time is concerned. To prove contraction in backward time it is
sufficient to prove it for points on \( W \) and then apply exactly
the same approximation argument. For \( z\in W\) this follows
immediately from the uniform expansivity assumption in forward time.
Indeed, letting \( v_{-n} = Df_{z}^{-n}(v) \), the expansivity
assumption gives 
\[ 
\|v\|\geq \|Df^{n}_{z_{-n}}(v_{-n})\|\geq
Ce^{\lambda n}\|v_{-n}\| 
\]
which immediately implies \( \|v_{-n}\|
\leq C^{-1}e^{-\lambda n}\|v\| \). 

It remains to show uniqueness of \( E^{u}(z) \) for each \( z\in
\Omega\). Suppose by contradiction that we could find two sequences 
\( z_{j}\to z \) and \( \tilde z_{j}\to z \) with corresponding
directions \( E^{u}_{z_{j}} \) and \( E^{u}_{\tilde z_{j}} \)
converging to two different directions \( E^{u}_{z} \) and \( \tilde
E^{u}_{z} \). Let \( v \in E^{u}_{z} \) and \( \tilde v \in \tilde
E^{u}_{z} \). Then \( v, \tilde v \) must be linearly independent and 
thus every other vector \( w \in T_{z}\mathbb R^{2} \) can be written 
as a linear combination \( w =  a_{1}v  + a_{2}\tilde v\) for some \(
w_{1}, w_{2} \in \mathbb R \). By linearity and the backward contraction estimates
obtained above this implies that 
\[
\|w_{-n}\| = \|Df^{-n}_{z}(w)\| \to 0
\]
as \( n\to \infty \). Since \( w \) was arbitrary this implies that
all vectors are shrinking to zero in backward time. But this is
impossible since we have assumed that \( |det Df| < 1 \) and thus \(
|det Df^{-1}| > 1 \)  on \( \Omega \). 
\end{proof}   

\begin{corollary}\label{unique}
    At every point \( z \in \Omega \) there exists a unique tangent space
    splitting \( E^{u}_{z}\oplus E^{s}_{z} \) which is invariant by
    the derivative \( Df \) and which satisfies the standard uniform
    hyperbolicity expansion and contraction estimates.
\end{corollary}
\begin{proof}
 Lemma \ref{limit} gives the expanding direction \( E^{u}_{z} \) of
 the splitting with the required hyperbolic expansion estimates in
 forward time. The invariance for points in \( W \) is automatic
 (since tangent directions to \( W \) are mapped to tangent directions
 to \( W \)), and the invariance for general points follows immediately
 from the definition of \( E^{u}_{z}= \lim E^{u}_{z_{j}} \), the
 invariance of \( E^{u}_{z_{j}} \) for \( z_{j}\in W \), and the
 continuity of \( Df \). 
 
 The stable direction \( E^{s}_{z} \) is given immediately by
as the limit of the sequence \( e^{(n)} \) of vectors most
 contracted by \( Df_{z}^{n} \), as discussed in section
 \ref{shadowing}. This also automatically gives the uniqueness and 
 invariance. 
 \end{proof}

To complete the proof of the Proposition, 
we just need to show that the given tangent
space splitting is continuous. This follows by standard arguments from 
the uniqueness proved in Corollary \ref{unique}. Indeed, for any
\( z\in\Omega \) and any sequence \( z_{j}\in\Omega \) with \(
z_{j}\to z \), every limit point of the corresponding sequence of
splittings \( E^{u}_{z_{j}}\oplus E^{s}_{z_{j}} \) must also be a
splitting \( \tilde E^{u}_{z}\oplus \tilde E^{s}_{z} \). By
approximation arguments identical to those used above it follows that 
this splitting must also satisfy the uniform hyperbolic contraction
and expansion estimates. Therefore. by uniqueness, it must coincide
with the existing splitting \( E^{u}_{z}\oplus E^{s}_{z} \). 
This completes the proof that \( \Omega \) is uniformly hyperbolic. 
 \end{proof}

\subsection{Lyapunov exponents for 
\protect\(  f_{a^{*}}\protect\)}
\label{Lyap}

\begin{figure}[h]
    \includegraphics[width=10cm]{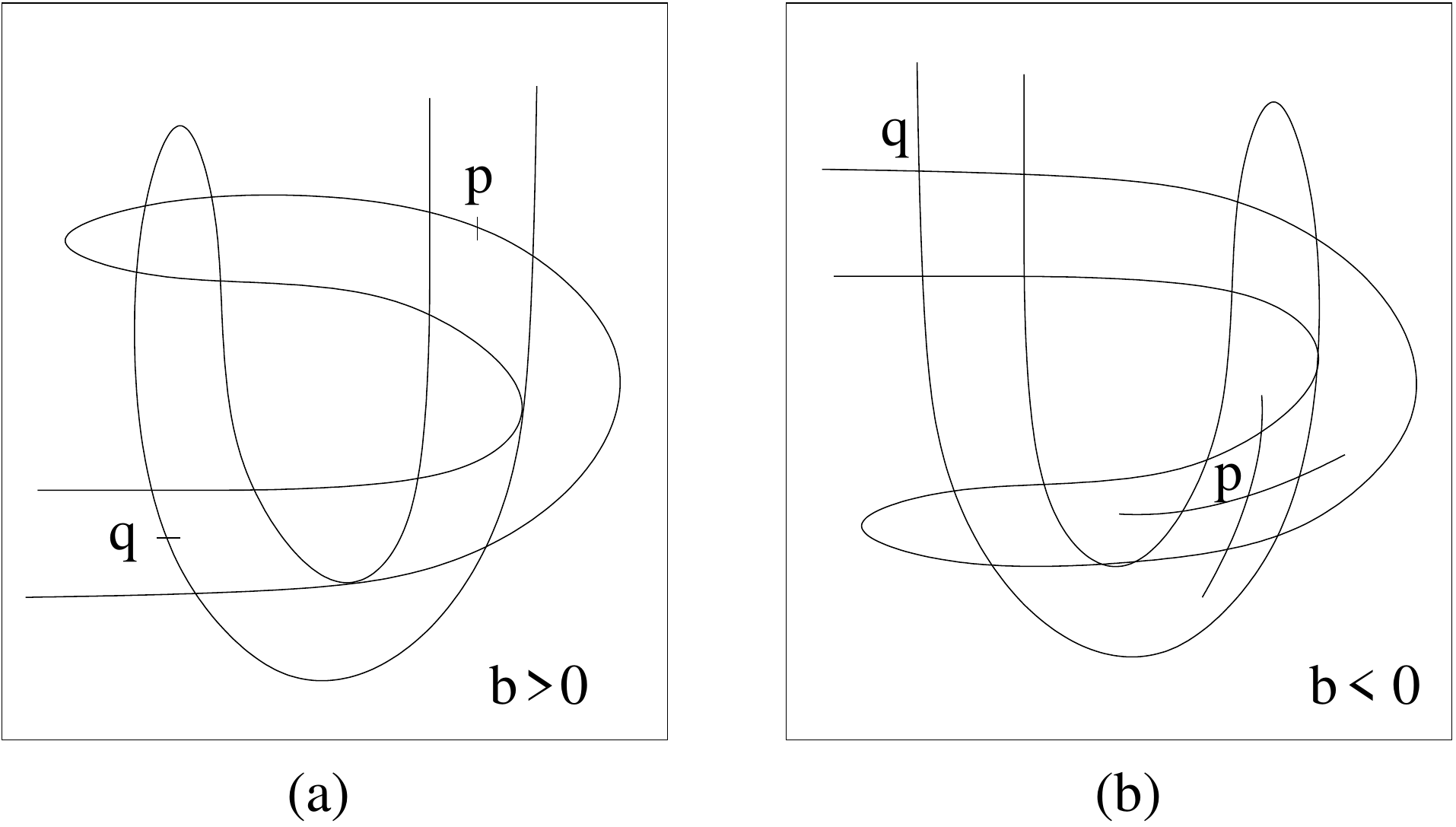}
 \caption{Invariant manifolds for \(a=a^{*}\)}\label{tangency}
 \end{figure}

Finally it remains to consider the dynamics of \( f_{a^{*}} \). 
Recall that \( a^{*} \) is defined on page \pageref{astar} 
as the first parameter
for which a tangency occurs between the compact parts of \( W^{s}(q) \)
and \( W^{u}(p) \), see Figure \ref{tangency} for the pictures in the 
two cases \( b>0 \) and \( b< 0 \). 

We need to show that, for
\( a=a^{*} \), all Lyapunov exponents are uniformly bounded away from 
0. 
We show that for each point
\( z \in \Omega_{a^{*}}\) \emph{not in the orbit of tangency } 
\( \mathcal T \)
(it is not necessary to consider the orbit of tangency since this is a
countable set without recurrence and can therefore not support any
invariant probability measure) 
there exists a constant \( C_{z} \), a
vector \( v_{z} \),  and a
sequence \( \{n_{i}\} \) 
with \( n_{i}\to\infty \) such that, for all \( i\geq 0 \),
\[
\|Df^{n_{i}}_{z}(v_{z})\|\geq C_{z}
e^{\lambda n_{i}}  \|v_{z}\|.
\]

This is obviously true
if the orbit of \( z \) never enters \( \Delta_{\varepsilon} \) in
forward  time or it enters \( \Delta_{\varepsilon} \) only a finite
number of times. Indeed suppose that there exists some \( k \) such
that \( f^{i}(z)
\notin\Delta_{\varepsilon} \) for all \( i\geq k \). Then let \( w \) 
be a vector which is mapped to the horizontal vector 
\( w_{k}= Df_{z}^{k}(w)  \) 
after \( k \) iterations.
Then by \eqref{UE1} we have 
\( \| Df_{z_{k}}^{k+n}(w) \| \geq C_{\varepsilon} e^{\lambda n}
\|w_{k}\|\) 
for all \( n\geq 1 \). This implies that there exists a constant \(
C_{z} \) such that 
\( \| Df_{z}^{k+n}(w) \| \geq C_{z} e^{\lambda (k+n)}
\|w\|\)  for all \( n\geq 1 \).

Otherwise there exists an infinite sequence \( 0< m_{1}< \cdots
< m_{k}< \cdots\) such that \( m_{k}\to\infty \) and \(
f^{m_{k}}(z)\in\Delta_{\varepsilon} \).  By Lemma \ref{unique
critical}, 
\( z_{m_{i}}=f^{m_{i}}(z) \) lies on either a \( C^{2} \) long admissible
curve or a \( C^{1+1} \) 
long admissible curve which is the \( C^{1} \) limit of \( C^{2} \)
long admissible curves in \( W^{u}(p) \). 
Since \( z \) has an infinite number of returns to \(
\Delta_{\varepsilon} \), this implies in particular that \(
z\notin W^{s}(q) \) and so \( z_{m_{i}}\notin W^{s}(q) \) and so 
 there exists \( n_{i}=n(z_{m_{i}}) \) such that \(
f(z_{m_{i}}) \in \mathcal V_{n_{i}}\setminus\mathcal V_{n_{i}+1} \). Therefore
\emph{exactly} the same arguments as in 
Lemmas \ref{binding} and \ref{limit} show that for a vector \( v_{i} \) tangent to
such an admissible curve \( \gamma \) at \( z_{m_{i}} \) we have 
\begin{equation}\label{lasteq} 
\|Df^{n_{i}+1}_{z_{m_{i}}}(v_{i})\|
\geq e^{\lambda (n_{i}+1)}\|v_{i}\|. 
\end{equation}

%

Notice that since the \( C^{1} \) limits of \( C^{2} \) admissible
curves are unique, as proved above, we have \(
v_{i+1}=Df^{m_{i+1}-m_{i}}(v_{i}) \). Then, by \eqref{UE1} 
and \eqref{lasteq} we have 
\[ 
\|Df^{m_{i}+n_{i}+1-m_{1}}(v_{1})\|\geq 
e^{\lambda(m_{i}+n_{i}+1-m_{1})}\|v_{1}\|.
\]
Then we can define \( v_{z}= Df^{-m_{1}}(v_{1}) \) 
and we have 
\(  \|Df^{m_{i}+n_{i}+1}(v_{z})\|\geq 
C_{z} e^{\lambda(m_{i}+n_{i}+1)}\|v_{z}\|\).
where the constant \( C_{z} \) is required simply to compensate for
the possible lack of expansion for the first \( n_{1} \) iterates. In 
particular it can be chosen by considering the maximum possible
contraction along the orbit of \( z \) for the first \( n_{1} \)
iterations
\[ 
C_{z} =  \min_{\|v\|=1} \|Df^{n_{1}}_{z}(v)\|.
\]
We have shown therefore that for each \( z\in\Omega \) 
\(
\limsup_{n\to\infty} \frac{1}{n}\ln \|Df^{n}_{z}\|\geq \lambda.
\)
This clearly implies the same bound for the limit wherever it exists. 
In particular any point which is typical for some ergodic invariant
probability measure and for which therefore such a limit does exist,
will have a positive Lyapunov exponent \( \geq \lambda \). By
dissipativity this immediately implies also that the other Lyapunov
exponent is negative and uniformly bounded away from 0 both in the
dynamical and in the parameter space.

  \end{document}